\documentclass[paper=a4,abstract=true,numbers=noenddot]{scrartcl}

\usepackage[english]{babel}  
\usepackage{hyperref}
\usepackage[intlimits]{amsmath}
\usepackage{amsthm}
\usepackage{amssymb} 
\usepackage{amsfonts,mathrsfs,bbm}
\allowdisplaybreaks[1]
\usepackage[numbers]{natbib}
\usepackage{graphicx}
\usepackage[space]{grffile} 
\usepackage{tikz}
\usepackage{pgfplots}
\usepackage[T1]{fontenc}
\usepackage{mathptmx}
\usepackage[scaled=.92]{helvet}
\usepackage{courier}

\usepackage{authblk}
\usepackage{amsmath,bm}
\usepackage{graphicx}
\usepackage{mathrsfs}
\usepackage{physics}
\usepackage{tensor}
\usepackage{xcolor}
\usepackage{enumerate}

\usepackage{subcaption}

\def\Xint#1{\mathchoice
   {\XXint\displaystyle\textstyle{#1}}%
   {\XXint\textstyle\scriptstyle{#1}}%
   {\XXint\scriptstyle\scriptscriptstyle{#1}}%
   {\XXint\scriptscriptstyle\scriptscriptstyle{#1}}%
   \!\int}
\def\XXint#1#2#3{{\setbox0=\hbox{$#1{#2#3}{\int}$}
     \vcenter{\hbox{$#2#3$}}\kern-.5\wd0}}
\def\ddashint{\Xint=}
\def\dashint{\Xint-}
\usetikzlibrary{patterns}
\pgfplotsset{compat=1.18}
\usepackage{siunitx}
\usepackage{array} 
\usepackage{booktabs} 
\usepackage{setspace}

\theoremstyle{remark}

\title{Partial integration based regularization in BEM for 3D elastostatic problems: The role of line integrals}

\author[1]{V. Lakshmi Keshava}
\author[1]{M. Schanz}
\affil[1]{Institute of Applied Mechanics, Graz University of
  Technology, Technikerstraße 4/II, 8010 Graz, Austria,
  v.lakshmikeshava@tugraz.at, m.schanz@tugraz.at}
\date{}                     
\begin{document}
	
\maketitle
\section*{Abstract}
The Boundary Element Method (BEM) is a powerful numerical approach for solving 3D elastostatic problems, particularly useful for crack propagation in fracture mechanics and half-space problems. A key challenge in BEM lies in handling singular integral kernels. Various analytical and numerical integration or regularization techniques address this, including one that combines partial integration with Stokes' theorem to reduce hyper-singular and strong singular kernels to weakly singular ones. This approach typically assumes a closed surface, omitting the boundary integrals from Stokes' theorem. In this paper, these usually neglected boundary line integrals are introduced and their significance is demonstrated, first in a pure half-space problem, and then shown to be redundant in fast multipole method (FMM) based BEM, where geometry partitioning produces pseudo open surfaces.

\textbf{Keywords}: Half-space problem; Fast multipole method; Partial integration; Stokes' theorem; Line integrals

\section{Introduction} \label{sec1}
The Boundary Element Method (BEM) is a widely-used numerical tool for solving linear elastostatic boundary value problems, particularly effective for half-space problems and crack-propogation problems in fracture mechanics. BEM offers certain advantages over the traditional Finite Element Method (FEM) in treating such complex geometries, mainly by reducing the dimensionality of the discretization and inherently fulfilling the radiation condition in unbounded domains. The method is well-established in both collocation and Galerkin approaches, and a detailed mathematical explanation can be found in \cite{brebbia84, steinbach08}. 

The primary challenge of BEM lies in the treatment of singular integral kernels that emerge due to the nature of the fundamental solution and its derivatives. In the collocation approach with standard boundary integral equations (SBIE), two types of singularities are encountered: weakly singular and strongly singular kernels. Weakly singular integrals are simply improper integrals that can be handled efficiently using coordinate transformations, which cancel out the singularities, such that standard Gaussian quadrature can be applied \cite{duffy82}. In contrast, strongly singular integrals should be defined with a limiting process, as a Cauchy princpal value, and cannot be computed directly using a standard numerical scheme. The symmetric Galerkin BEM or the collocation BEM with hypersingular boundary integral equations (HBIE)  introduces an additional hypersingular integral, which must be evaluated as finite part integral in the sense of Hadamard \cite{hadamard1923}. Like the strongly singular kernels, hypersingular integrals also poses challenges within the standard numerical scheme.   

The proper treatment of these singular kernels is well understood in the BEM community, and several techniques have been developed to address this challenge. These techniques can be broadly classified into analytical and numerical regularization, depending on how the singularity is handled. A comprehensive review of the different regularization methods for boundary value problems in continuum mechanics can be found in \cite{tanaka94}. 

Analytical regularization evaluates singular integrals exactly, providing high accuracy. However, its implementation is often complex, requiring intricate formulae that depend on the problem’s geometry and the choice of trial functions, thereby limiting its general applicability. One of the earliest analytical treatments of singular integrals in collocation BEM was introduced by Guiggiani \textit{et al}. \cite{guiggiani90, guiggiani91, guiggiani98}. Their approach involved temporarily deforming the boundary around the singularity and translating the singular integrands to intrinsic coordinates, allowing for an exact limiting process. This method, though effective, requires a Taylor series expansion of the kernel function, necessitating special treatment for each specific integral and introducing truncation errors. Gao \cite{gao10} introduced a different technique in which the non-singular part is expressed as a power series expansion, while the singularities are eliminated using the radial integration method \cite{gao02}. This approach is restricted to singular integrals in the Cauchy principal value sense and is formulated specifically for quadrilateral boundary elements. More recently, a kernel separation method was proposed by Xie \textit{et al}. \cite{Xie20} where the integrals are separated into two parts, which include the singular part and the regular part, where the singular part is computed through the Taylor expansion around the singular point. For nearly-singular integrals, a variety of analytical and semi-analytical approaches have been proposed; see, e.g., \cite{Zhong24, Han22}.

For the Galerkin approach, Hackbusch \& Sauter \cite{hackbush93} introduced an analytical integration technique using special variable substitutions and coordinate transformations to first regularize the integrals to weak singularity before applying the Duffy transformation. Andrä \& Schnack \cite{andra97} extended these formulae to curved elements by introducing new kernel functions with reduced singularity order, obtained via integration by parts of the standard Kelvin solution. Haas \& Kuhn \cite{haas02} further refined this approach for curved elements, employing a decomposed Kelvin fundamental solution without introducing additional kernel functions. In addition to these methods, numerous other analytical techniques have been developed, each offering distinct advantages depending on the problem at hand. Nevertheless, they all involve intricate mathematical formulae, making their implementation challenging. Most importantly, all the analytical formulae are derived in the intrinsic coordinate space (post-discretization), making them specific to the type of boundary elements used.

On the other hand, purely numerical techniques for computing the strong- and hypersingular integrals, while easier to implement, are relatively uncommon. The few existing formulae for Hadamard finite part integrals rely on obscure or insufficient parameters to ensure the stability of the numerical scheme, which can limit their practical effectiveness \cite{kutt75, mantic94}. A more practical and widely used approach is the hybrid approach, which combines classical Stokes' theorem with partial integration. In this approach, the tangential derivatives are 'shifted' from the decomposed fundamental solution introduced by Kupradze \textit{et al}. \cite{kupradze79} to the Cauchy data and/or the test function, making the integral weakly singular. The Duffy coordinate transformation \cite{duffy82} and standard Gauss quadrature are then applied to evaluate the integrals.

Although, Kupradze \textit{et al}. only dealt with strongly singular integrals, their method was later extended to hypersingular integrals by Sl\'{a}dek and Sl\'{a}dek \cite{sladek82}. Nedelec \cite{nedelec82} derived the corresponding regularization for the Galerkin formulation, which was further simplified by Han \cite{han94} by restricting the formulae to the isotropic case.

In this work, the formulae by Kupradze \textit{et al}. and Han are directly incorporated to handle singular integrals, primarily due to their ease of implementation. An additional advantage of this approach is that, because the integrals are decomposed and regularized, even for the regular case, the number of required integration points is significantly reduced (see, eg., \cite{kielhorn09}) compared to other analytical techniques, where regular integrals are left untreated. Moreover, nearly-singular integrals become less critical, and can be computed without additional treatment. However, this work does not claim that the hybrid method based on Stokes' theorem and integration by parts is inherently superior to other available techniques, as such comparisons are beyond its scope.

In the regularization process using Stokes’ theorem, the line integrals on the right-hand side are typically omitted under the assumption of a closed surface. However, this assumption fails for open surfaces, such as in half-space problems. In this paper, these omitted line integrals are explicitly presented for both strong and hypersingular cases, and their significance is demonstrated in the context of half-space problems. An alternative approach to address this issue was introduced in Kielhorn’s thesis \cite{kielhorn09}, where, instead of computing the line integrals, he proposed the use of so-called infinite elements in the half-space.

Another major challenge of BEM is its quadratic complexity in both storage and computation time. Several fast methods have been introduced to address this issue, beginning with Rokhlin’s work \cite{rokhlin85}, which reduced the complexity to logarithmic. Building on this, Greengard \& Rokhlin \cite{greengard87} developed the Fast Multipole Method (FMM), which has since been widely adopted in the boundary element framework. Liu \cite{liu09} published a comprehensive textbook on FMM with applications in collocation BEM, while Of et al. \cite{of05} applied it to Galerkin BEM in elastostatics. A more detailed literature review can be found in \cite{nishimura02}. In this work, a Chebyshev-interpolation-based FMM, originally introduced by Fong \& Darve \cite{fong09} as a black-box method, is used with minimal modifications from its original formulation.

The omission of line integrals in the regularization process also raises concerns in the FMM framework, where geometric clustering divides the domain into near- and far-field regions, seemingly violating the closed-surface assumption. In this paper, the omitted line integrals are incorporated using two different approaches and compared with the standard formulation without them. Numerical examples ultimately demonstrate that the line integrals are unnecessary in the context of FMM.

In Section \ref{sec2}, the basic governing equations are recalled and the corresponding collocation and Galerkin BEM formulations are presented. The regularization process here is only an extension of the work by Kupradze \textit{et al} \cite{kupradze79} and Han \cite{han94}, and as such, only the corresponding line integrals of the regularized double layer potential and the hypersingular operator are described in Section \ref{sec3}. Two versions of FMM are introduced in Section \ref{sec4}: first, the standard FMM with minimal modifications from \cite{fong09} and the second, the regularized FMM presented in  this paper, to address the problem of the near-field being an open-surface. To test the significance of the proposed line integrals and the regularized FMM, a few example problems are considerd in Section \ref{sec5}.  

\section{Boundary element formulation in elastostatics}
\label{sec2}

\subsection{Problem setting}\label{PS}
Consider a bounded Lipschitz domain $\Omega \subset \mathbb{R}^{3}$ with its boundary denoted by $\Gamma = \partial\Omega$. The linear isotropic elastostatics problem with no body forces is considered, with the displacement field denoted by $\mathbf{u}\left(\mathbf{x}\right)$. The governing equation for the considered boundary value problem can be written as 
\begin{equation}\label{goveq}
\left( \lambda + \mu \right)\nabla \left( \nabla \cdot\mathbf{u}\left ( \mathbf{x} \right ) \right) + \mu \nabla^2\mathbf{u}\left ( \mathbf{x} \right ) = 0\; \; \;\; \; \; \; \; \;\; \; \textup{for}  \: \: \mathbf{x}\in\Omega
\end{equation}
with $\lambda$ and $\mu$ being the Lam\'e parameters which are related to the Young's modulus $E$ and the Poisson's ratio $\nu$ of the material
\begin{equation}\label{Lame}
\lambda = \frac{E\nu}{(1+\nu)(1-2\nu)}, \quad \mu=\frac{E}{2(1+\nu)}.
\end{equation}
The $\nabla$-operator is used here to denote the spatial derivatives, with $\nabla=\textup{grad}$ and $\nabla \cdot=\textup{div}$. Considering a mixed boundary value problem, the Dirichlet and the Neumann boundary conditions are given by 
\begin{equation}\label{BCs}
\begin{aligned}
\mathbf{u}\left ( \mathbf{x} \right ) &= \mathbf{g}_{D}\left ( \mathbf{x} \right )\; \; \; \textup{for}  \: \: \mathbf{x}\in\Gamma_{D} \\
\mathbf{t}\left ( \mathbf{x} \right ) = \mathcal{T}_{\mathbf{x}}\left ( \mathbf{u}\left ( \mathbf{x} \right ) \right )&=\mathbf{g}_{N}\left ( \mathbf{x} \right )\; \; \; \textup{for}  \: \: \mathbf{x}\in\Gamma_{N}.
\end{aligned}
\end{equation}
The Lipschitz boundary $\Gamma$ is split into mutually disjoint sets $\Gamma_D$ and $\Gamma_N$ such that $\Gamma = \Gamma_D \cup \Gamma_N$, where the subscript $D$ and $N$ denoting the Dirichlet and the Neumann boundaries, respectively. The traction vector is represented by $\mathbf{t}\left ( \mathbf{x} \right )$, while $\mathcal{T}_{\mathbf{x}}$ is the elastic traction operator defined as 
\begin{equation}
\mathcal{T}_{\mathbf{x}}\left ( \mathbf{u}\left ( \mathbf{x} \right ) \right ):=  \boldsymbol{\sigma}\left ( \mathbf{u},\mathbf{x} \right )\mathbf{n}\left ( \mathbf{x} \right ),
\end{equation}
with the stress tensor $\boldsymbol{\sigma}\left(\mathbf{u},\mathbf{x}\right)$ and the outward normal unit vector $\mathbf{n}\left ( \mathbf{x} \right )$.

\subsection{Boundry integral equation}
Both the indirect and direct formulations are possible, but only the latter is considered here. The fundamental solution for \eqref{goveq} is given by the Kelvin tensor \cite{han94} 
\begin{equation}\label{FS}
\mathbf{U}_{ij}\left ( \mathbf{x},\mathbf{y} \right )=\frac{1}{8\pi\mu(\lambda + 2 \mu)}\left [ \frac{\left (\lambda + 3 \mu \right )}{\left | \mathbf{x}-\mathbf{y} \right |} \delta_{ij}+(\lambda + \mu)\frac{\left ( x_{i}-y_{i} \right )\left ( x_{j}-y_{j} \right )}{\left | \mathbf{x}-\mathbf{y} \right |^{3}}\right ]
\end{equation}
for $i,j=1,2,3$, where  $\delta_{ij}$ denotes the Kronecker delta. Based on Betti's second formula, the standard boundary integral equation (SBIE) is given by 
\begin{equation}\label{SBIE}
\mathbf{C(x)}\mathbf{u(x)}= \int_\Gamma\mathbf{U(x,y)}\;\mathbf{t(y)}\mathrm{d}s_\mathbf{y} - \int_{\Gamma} (\mathcal{T}_{\mathbf{y}}\mathbf{U(x,y)})^\top\mathbf{u(y)} \mathrm{d}s_\mathbf{y} \;\;\;\;\;\;\;  \forall \: \mathbf{x,y}\in \Gamma
\end{equation}
where $\mathbf{C(x)}$ is the so called integral-free term, which is dependant on the Poisson's ratio $\nu$ (in elasticity problems) and the geometry (the solid angle) at $\mathbf{x}$, in collocation BEM \cite{gaul03, mantic93}. When the boundary $\Gamma$ is at least differentiable at $\mathbf{x}$, and in the Galerkin formulation, this term simplifies to $\mathbf{C(x)} = \frac{1}{2}\mathbf{\mathbf{I}}$, where $\mathbf{I}$ represents the identity matrix. The hypersingular boundary integral equation (HBIE) is written as 
\begin{equation}\label{HBIE}
\begin{aligned}
\left[\mathbf{I}-\mathbf{C(x)}\right]\mathbf{t(x)}&= \mathcal{T}_{\mathbf{x}}\int_\Gamma\mathbf{U(x,y)\;t(y)}\mathrm{d}s_\mathbf{y} \\
&- \mathcal{T}_{\mathbf{x}}\int_\Gamma (\mathcal{T}_{\mathbf{y}}\mathbf{U(x,y)})^\top\mathbf{u(y)} \mathrm{d}s_\mathbf{y} \;\;\;\;\;\;\;  \forall \: \mathbf{x,y}\in \Gamma.
\end{aligned}
\end{equation}
Corresponding to the above SBIE and the HBIE, the standard boundary integral opearators are introduced. The single layer potential is defined as
\begin{equation}\label{SLP}
(\mathcal{V}\mathbf{t})(\mathbf{x}):= \int_\Gamma\mathbf{U(x,y)}\;\mathbf{t(y)}\mathrm{d}s_\mathbf{y}.
\end{equation}
The double layer potential and the adjoint double layer potential are given by
\begin{align}\label{DLP}
(\mathcal{K}\mathbf{u})(\mathbf{x})&:= \int_\Gamma (\mathcal{T}_{\mathbf{y}}\mathbf{U(x,y)})^\top\mathbf{u(y)} \mathrm{d}s_\mathbf{y},\\
(\mathcal{K'}\mathbf{t})(\mathbf{x})&:= \mathcal{T}_{\mathbf{x}}\int_\Gamma\mathbf{U(x,y)\;t(y)}\mathrm{d}s_\mathbf{y}.
\label{ADLP}
\end{align}
Finally, the hypersingular operator is defined as
\begin{equation}\label{HSO}
(\mathcal{D}\mathbf{u})(\mathbf{x}):= - \mathcal{T}_{\mathbf{x}}\int_\Gamma (\mathcal{T}_{\mathbf{y}}\mathbf{U(x,y)})^\top\mathbf{u(y)} \mathrm{d}s_\mathbf{y}.
\end{equation}
As the name suggests, the hypersingular operator \eqref{HSO} exhibits hypersingular behavior as $\mathbf{x}\to\mathbf{y}$ and must, therefore, be interpreted as a finite part integral in the sense of Hadamard \cite{hadamard1923}. Similarly, the double layer potential \eqref{DLP} and its adjoint \eqref{ADLP} exhibit strong singularities as $\mathbf{x}\to\mathbf{y}$ and must be understood in the sense of a Cauchy principal value. In the literature, the symbols $\dashint$ and $\ddashint$ are often used to denote the Cauchy principal value and the finite-part integral in the sense of Hadamard, respectively. However, for simplicity, such notation is omitted here.

\subsection{Spatial discretization}
The boundary $\Gamma$ is discretized using flat linear triangles, resulting in an approximate boundary $\Gamma^h$, expressed as
\begin{equation}\label{tri}
\Gamma \approx \Gamma^h = \bigcup_{k=1}^{K}\tau_k.
\end{equation}
Thus, $\Gamma^h$ consists of a union of $K$ triangular elements $\tau_k$. The unknown fields $\mathbf{u}$ and $\mathbf{t}$ are then approximated using piecewise linear continuous basis functions $\varphi_j^1$ and piecewise constant basis functions $\varphi_j^0$ respectively, as
\begin{equation}\label{discretization}
\mathbf{u(x)} \approx \sum_{j=1}^M \mathbf{u}_j \varphi_j^1(\mathbf{x})\;\;\;\;\;\;\;\;\;\textup{and}\;\;\;\;\;\;\;\;\; \\
\mathbf{t(x)} \approx \sum_{j=1}^N \mathbf{t}_j \varphi_j^0(\mathbf{x}).
\end{equation}
This leads to $M$ Dirichlet unknowns $\mathbf{u}_j$ and $N$ Neumann unknowns $\mathbf{t}_j$. Substituting these approximations into the boundary integral equations leads to a system of semi-discrete equations, which can be solved using either the collocation or Galerkin approach. Since both methods are employed in this study to examine their behavior in half-space problems, they are briefly outlined in the following section.

\subsection{Boundary element formulations}\label{sec2.4}
The mixed boundary conditions \eqref{BCs} imply that the unknowns $\mathbf{u(x)}$ and $\mathbf{t(x)}$ are partially known, depending on the boundary segment. To apply the integral operators on the known and unknown Cauchy data separately, they are decomposed as
\begin{equation}
\mathbf{u(x)} = \mathbf{\tilde{u}(x)} + \mathbf{\tilde{g}}_D(\mathbf{x})\;\;\;\;\;\textup{and}\;\;\;\;\; \\
\mathbf{t(x)} = \mathbf{\tilde{t}(x)} + \mathbf{\tilde{g}}_N(\mathbf{x})\;\;\;\;\; \textup{for } \mathbf{x}\in \Gamma,
\end{equation}
where $\mathbf{\tilde{u}(x)}$ and $\mathbf{\tilde{t}(x)}$ are the unknown Dirichlet and Neumann data, respectively. The prescribed Cauchy data $\mathbf{g}_D(\mathbf{x})$ and $\mathbf{g}_N(\mathbf{x})$ are extended to the arbitrary fixed functions $\mathbf{\tilde{g}}_D(\mathbf{x})$ and $\mathbf{\tilde{g}}_N(\mathbf{x})$, chosen to vanish on the complementary boundary segments \citep{steinbach08}. 

With this decomposition, and substituting the spatial shape functions \eqref{discretization} into the SBIE \eqref{SBIE}, the discretized equation can be written in matrix form as
\begin{equation}\label{colloMatrix}
\begin{bmatrix}
 \mathsf{V}_D & -\mathsf{K}_D \\
 \mathsf{V}_N & -\mathsf{K}_N
\end{bmatrix} \begin{bmatrix}
\mathbf{\tilde{t}} \\
\mathbf{\tilde{u}}
\end{bmatrix} = \begin{bmatrix}
\mathsf{f}_D \\
\mathsf{f}_N
\end{bmatrix},
\end{equation} 
where $\mathsf{V}$ and $\mathsf{K}$ denote the discrete single-layer and double-layer collocation matrices, respectively, given by 
\begin{equation*}
\begin{aligned}
\mathsf{V}_D[i,j] &= \int_{supp(\varphi_j^0)}\mathbf{U}(\mathbf{x}_i\mathbf{,y})\;\varphi_j^0(\mathbf{y})\;\mathrm{d}s_{\mathbf{y}}\;\;\;\;\;\;\;\;\;\;\;\;\;\;\;\;\;\;\;\;\;\;\;\;\;\;\textup{for }\mathbf{x}_i\in \Gamma_D \\
\mathsf{V}_N[i,j] &= \int_{supp(\varphi_j^0)}\mathbf{U}(\mathbf{x}_i\mathbf{,y})\;\varphi_j^0(\mathbf{y})\;\mathrm{d}s_{\mathbf{y}}\;\;\;\;\;\;\;\;\;\;\;\;\;\;\;\;\;\;\;\;\;\;\;\;\;\;\textup{for }\mathbf{x}_i\in \Gamma_N \\
\mathsf{K}_D[i,j] &= \mathbf{C}(\mathbf{x}_i) +  \int_{supp(\varphi_j^1)}(\mathcal{T}_{\mathbf{y}}\mathbf{U}(\mathbf{x}_i\mathbf{,y}))^\top\;\varphi_j^1(\mathbf{y})\;\mathrm{d}s_{\mathbf{y}}\;\;\;\;\;\textup{for }\mathbf{x}_i\in \Gamma_D \\
\mathsf{K}_N[i,j] &= \mathbf{C}(\mathbf{x}_i) +  \int_{supp(\varphi_j^1)}(\mathcal{T}_{\mathbf{y}}\mathbf{U}(\mathbf{x}_i\mathbf{,y}))^\top\;\varphi_j^1(\mathbf{y})\;\mathrm{d}s_{\mathbf{y}}\;\;\;\;\;\textup{for }\mathbf{x}_i\in \Gamma_N.
\end{aligned}
\end{equation*}
The right hand side vectors $\mathsf{f}_D$ and $\mathsf{f}_N$ are computed analogously. The collocation points $\mathbf{x}_i$ are chosen as the nodes of the linear triangular elements for $\mathbf{x}_i\in \Gamma_N$ and as the midpoints of the constant triangular elements for $\mathbf{x}_i\in \Gamma_D$.

For the symmetric Galerkin method, both the SBIE \eqref{SBIE} and the HBIE \eqref{HBIE} are tested with appropriate test functions $\mathbf{v(x)}$ and $\mathbf{w(x)}$, leading to the weak formulation
\begin{equation}\label{VarForm}
\begin{aligned}
\langle \mathcal{V}\mathbf{\tilde{t}}, \mathbf{v} \rangle_{\Gamma_D} - \langle \mathcal{K}\mathbf{\tilde{u}}, \mathbf{v} \rangle_{\Gamma_D} &= \langle \tfrac{1}{2} \hat{I}\mathbf{\tilde{g}}_D+\mathcal{K} \mathbf{\tilde{g}}_D - \mathcal{V} \mathbf{\tilde{g}}_N, \mathbf{v} \rangle_{\Gamma_D} \\
\langle \mathcal{K'}\mathbf{\tilde{t}}, \mathbf{w} \rangle_{\Gamma_N} - \langle \mathcal{D}\mathbf{\tilde{u}}, \mathbf{w} \rangle_{\Gamma_N} &= \langle \tfrac{1}{2} \hat{I}\mathbf{\tilde{g}}_N-\mathcal{K'}  \mathbf{\tilde{g}}_N - \mathcal{D} \mathbf{\tilde{g}}_D, \mathbf{w} \rangle_{\Gamma_N},
\end{aligned}
\end{equation}
where the known and unknown Cauchy data are separated as in the collocation approach. The operator $\hat{I}$ denotes the identity operator 
\begin{equation}\label{IOperator}
(\hat{I}w)(\mathbf{x}) := \int_{\Gamma}\delta(\mathbf{y-x})w(\mathbf{y})\mathrm{d}s_{\mathbf{y}}. 
\end{equation}
The inner product $$\langle \mathbf{f},\mathbf{g}\rangle_\Gamma:=\int_\Gamma \mathbf{f(x)g(x)} \: \mathrm{d}\mathbf{x}$$ implies that the Galerkin formulation requires double integrations. Inserting the spatial approximations \eqref{discretization} into the above variational formulation yields the linear system
\begin{equation}\label{GalMatrix}
\begin{bmatrix}
 \mathsf{\widehat{V}} & -\mathsf{\widehat{K}} \\
 \mathsf{\widehat{K}}^\top & \mathsf{\widehat{D}}
\end{bmatrix} \begin{bmatrix}
\mathbf{\tilde{t}} \\
\mathbf{\tilde{u}}
\end{bmatrix} = \begin{bmatrix}
\mathsf{\widehat{f}}_D \\
\mathsf{\widehat{f}}_N
\end{bmatrix},
\end{equation}
where $\mathsf{\widehat{V}}$ denotes the disctrete single layer Galerkin matrix given by
\begin{equation}\label{Vgal}
\mathsf{\widehat{V}}[i,j] = \int_{supp(\varphi_i^0)}\varphi_i^0(\mathbf{x}) \;\int_{supp(\varphi_j^0)}\mathbf{U}(\mathbf{x}\mathbf{,y})\;\varphi_j^0(\mathbf{y})\;\mathrm{d}s_{\mathbf{y}}\;\mathrm{d}s_{\mathbf{x}}.
\end{equation} 
$\mathsf{\widehat{K}}$ and $\mathsf{\widehat{K}}^\top$ are the discrete double layer and discrete adjoint double layer Galerkin matrices
\begin{equation}\label{Kgal}
\begin{aligned}
\mathsf{\widehat{K}}[i,j] &= \int_{supp(\varphi_i^0)}\varphi_i^0(\mathbf{x}) \;\int_{supp(\varphi_j^1)}(\mathcal{T}_{\mathbf{y}}\mathbf{U}(\mathbf{x}\mathbf{,y}))^\top\;\varphi_j^1(\mathbf{y})\;\mathrm{d}s_{\mathbf{y}}\;\mathrm{d}s_{\mathbf{x}}, \\
\mathsf{\widehat{K}}^\top[i,j] &= \int_{supp(\varphi_i^1)}\varphi_i^1(\mathbf{x}) \mathcal{T}_{\mathbf{x}}\;\int_{supp(\varphi_j^0)}\mathbf{U}(\mathbf{x}\mathbf{,y})\;\varphi_j^0(\mathbf{y})\;\mathrm{d}s_{\mathbf{y}}\;\mathrm{d}s_{\mathbf{x}}, 
\end{aligned}
\end{equation}
and $\mathsf{\widehat{D}}$ is the discrete hypersingular operator matrix given by
\begin{equation}\label{Dgal}
\mathsf{\widehat{D}}[i,j] = - \int_{supp(\varphi_i^1)}\varphi_i^1(\mathbf{x}) \mathcal{T}_{\mathbf{x}}\;\int_{supp(\varphi_j^1)}(\mathcal{T}_{\mathbf{y}}\mathbf{U}(\mathbf{x}\mathbf{,y}))^\top\;\varphi_j^1(\mathbf{y})\;\mathrm{d}s_{\mathbf{y}}\;\mathrm{d}s_{\mathbf{x}}.
\end{equation}
The right-hand side vectors are computed analogously, similar to the collocation method. Weakly singular integrals in the Galerkin formulation are handled numerically using the approach by Erichsen and Sauter \cite{erichsen98}, while those in the collocation method are treated via the Duffy transformation \citep{duffy82}. Strongly singular and hypersingular integrals are first transformed into weakly singular ones using partial integration and Stokes' theorem before numerical evaluation.

\section{Partial integration based regularization}
\label{sec3}

The regularization procedure adopted here follows the works of Han \cite{han94}. As previously mentioned, Han’s work is extended to accommodate cases where the line integrals arising from Stokes’ theorem are retained rather than omitted. In this section, the regularized kernels for the strongly singular double-layer potential and the hypersingular operator are presented, explicitly incorporating these line integrals. 

\subsection{The double layer potential $\mathcal{K}$}
Applying Stokes' theorem and integrating by parts, the double layer potential of linear elastostatics can be reformulated as
\begin{equation}\label{newDLP}
\begin{aligned}
 \left<\mathcal{K}\mathbf{\tilde{u},v} \right> &= 2\mu \int_{\Gamma}\mathbf{v\left ( x \right )}\cdot\int_{\Gamma}\mathbf{U}\cdot\left ( \mathcal{M}_\mathbf{y} \cdot\mathbf{\tilde{u}\left ( y \right )} \right )\mathrm{d}s_{\mathbf{y}} \mathrm{d}s_{\mathbf{x}} \\
 &- \frac{1}{4 \pi}\int_{\Gamma}\mathbf{v\left ( x \right )}\cdot\int_{\Gamma}\frac{1}{\left | \mathbf{x}-\mathbf{y} \right |}\left ( \mathcal{M}_\mathbf{y}   \cdot\mathbf{\tilde{u}\left (y  \right )}\right )\mathrm{d}s_{\mathbf{y}}\mathrm{d}s_{\mathbf{x}} \\
&+\frac{1}{4 \pi}\int_{\Gamma}\mathbf{v\left ( x \right )} \cdot\int_{\Gamma} \frac{\partial}{\partial \mathbf{n\left ( y \right )}}\frac{1}{\left | \mathbf{x}-\mathbf{y} \right |} \mathbf{\tilde{u}\left ( y \right )}\mathrm{d}s_{\mathbf{y}}\mathrm{d}s_{\mathbf{x}} \\
 & - \frac{1}{4 \pi} \int\nolimits_{\Gamma}\mathbf{v\left ( x \right )}\cdot\int\nolimits_{\partial \Gamma}\frac{1}{\left | \mathbf{x}-\mathbf{y} \right |} \left ( \mathbf{\tilde{u}}\left ( \mathbf{y} \right ) \times \mathrm{d}\boldsymbol\zeta_{\mathbf{y}}\right )\mathrm{d}s_{\mathbf{x}} \\
      &+ 2\mu \int\nolimits_{\Gamma}\mathbf{v\left ( x \right )}\cdot \int\nolimits_{\partial \Gamma}\left[\mathbf{U}\circledast \mathbf{\tilde{u}}\left ( \mathbf{y} \right )\right]^\top \cdot \mathrm{d}\boldsymbol\zeta_{\mathbf{y}}\;\mathrm{d}s_{\mathbf{x}},
\end{aligned}
\end{equation} 
where $\mathcal{M}$ is the Günter derivative operator, defined as 
\begin{equation*}
\mathcal{M}_{\mathbf{y}}:=\begin{pmatrix}
0 & n_2 \frac{\partial }{\partial y_1} - n_1\frac{\partial }{\partial y_2} & n_3 \frac{\partial }{\partial y_1} - n_1\frac{\partial }{\partial y_3} \\
n_1 \frac{\partial }{\partial y_2} - n_2\frac{\partial }{\partial y_1} & 0  & n_3 \frac{\partial }{\partial y_2} - n_2\frac{\partial }{\partial y_3} \\
 n_1 \frac{\partial }{\partial y_3} - n_3\frac{\partial }{\partial y_1}& n_2 \frac{\partial }{\partial y_3} - n_3\frac{\partial }{\partial y_2} & 0
\end{pmatrix},
\end{equation*}
and $$\mathbf{U} \circledast \mathbf{\tilde{u}}\left ( \mathbf{y} \right ) := \left[\mathbf{U}_{;,1}\times \mathbf{\tilde{u}}\left ( \mathbf{y} \right ) \quad \mathbf{U}_{;,2}\times \mathbf{\tilde{u}}\left ( \mathbf{y} \right )\quad \mathbf{U}_{;,3}\times \mathbf{\tilde{u}}\left ( \mathbf{y} \right ) \right].$$ 
The symbols $(\cdot)$ and $(\times)$ denote the standard vector dot product and cross product, respectively. The notation $\partial\Gamma$ denotes the boundary of the surface $\Gamma$, and the integral over $\partial\Gamma$ is therefore a line integral. The differential term $\mathrm{d}\boldsymbol\zeta_{\mathbf{y}}$ in the line integrals is indeed a vector quantity. This formulation expresses the double layer potential as a combination of the single and double layer potentials from the Laplace equation, along with the single layer potential in linear elastostatics. Additionally, two line integral terms arise from the right hand side of the Stokes' theorem. 

All differential operators in the first three surface integral terms are shifted onto $\mathbf{\tilde{u}\left ( y \right )}$ by reformulating the traction operator and applying Stokes' theorem and partial integration. The more widely recognized form of the regularized double layer potential \cite{kupradze79, han94} omits the two line integrals, as these terms vanish for closed surfaces. The regularized adjoint double layer potential can be written in a similar manner, and its representation is omitted here for brevity.

\subsection{The hypersingular operator $\mathcal{D}$}

Analogous to the double-layer potential, the regularized hypersingular operator in linear elastostatics can be expressed as
\begin{equation}\label{newHSO}
\begin{aligned}
\left<\mathcal{D}\mathbf{\tilde{u},v} \right> &= \frac{2\mu}{4 \pi} \int\nolimits_{\Gamma}\int\nolimits_{\Gamma}\frac{1}{\left | \mathbf{x}-\mathbf{y} \right |} \sum_{k,i=1}^{3} \frac{\partial v_i\left ( \mathbf{x} \right ) } {\partial S_k \left ( \mathbf{x} \right )} \; \frac{\partial \tilde{u}_i\left ( \mathbf{y} \right ) }{\partial S_k \left ( \mathbf{y}\right )} \;\;\mathrm{d}s_{\mathbf{y}} \mathrm{d}s_{\mathbf{x}} \\ 
 &- \frac{\mu}{4 \pi} \int\nolimits_{\Gamma}\int\nolimits_{\Gamma}\frac{1}{\left | \mathbf{x}-\mathbf{y} \right |} \; \; \frac{\partial }{\partial \mathbf{S} \left ( \mathbf{x} \right )}\cdot \mathbf{v\left ( x \right )} \;\; \frac{\partial  }{\partial \mathbf{S} \left ( \mathbf{y} \right )} \cdot \mathbf{\tilde{u}\left ( y \right )}\;\;\mathrm{d}s_{\mathbf{y}}\mathrm{d}s_{\mathbf{x}} \\
 & - \frac{\mu}{4 \pi} \int\nolimits_{\Gamma} \int\nolimits_{\Gamma} \left ( \mathcal{M}_\mathbf{x}\cdot\mathbf{v\left ( x \right )} \right )\cdot \left ( 4\mu \mathbf{U}-\frac{2}{\left | \mathbf{x}-\mathbf{y} \right |} \mathbf{I}_3 \right ) \cdot   \left ( \mathcal{M}_\mathbf{y} \cdot\mathbf{\tilde{u}\left ( y \right )} \right ) \;\;\mathrm{d}s_{\mathbf{y}}\mathrm{d}s_{\mathbf{x}}  \\
&- \frac{\mu}{4 \pi} \int\nolimits_{\Gamma}\int\nolimits_{\partial \Gamma} \frac{1}{\left | \mathbf{x}-\mathbf{y} \right |} \left[\left(\frac{\partial }{\partial \mathbf{S} \left ( \mathbf{y} \right)} \otimes \mathbf{\tilde{u}\left ( y \right )}\right)\cdot\mathbf{v\left ( x \right )}\right]\cdot \mathrm{d}\boldsymbol\zeta_{\mathbf{x}}\;\;\mathrm{d}s_{\mathbf{y}}\\
&+\frac{\mu}{4 \pi} \sum_{l=1}^{3}\sum_{i,j=1}^{3}\int\nolimits_{\Gamma}\left[ \mathcal{M}_{\mathbf{x}_{li}}  v_j \left(\mathbf{x}\right)\right]\int\nolimits_{\partial\Gamma}\sum_{k=1}^{3}\epsilon_{ijk}\tilde{u}_l\left(\mathbf{y}\right)\frac{1}{\left | \mathbf{x}-\mathbf{y} \right |}\mathbf{e}_k \cdot \mathrm{d}\boldsymbol\zeta_{\mathbf{y}} \;\; \mathrm{d}s_{\mathbf{x}} \\
&- \frac{\mu}{4 \pi} \int\nolimits_{\Gamma}\mathbf{\tilde{u}}\left(\mathbf{y}\right)\cdot\int\nolimits_{\partial\Gamma}\sum_{j=1}^{3}\begin{bmatrix}
\left(\mathcal{M}_{\mathbf{x}_{3j}}\frac{1}{\left | \mathbf{x}-\mathbf{y} \right |} \; v_j \left(\mathbf{x}\right) \; \mathbf{e}_2 -\mathcal{M}_{\mathbf{x}_{2j}}\frac{1}{\left | \mathbf{x}-\mathbf{y} \right |} \; v_j \left(\mathbf{x}\right) \; \mathbf{e}_3\right)\cdot \mathrm{d}\boldsymbol\zeta_{\mathbf{x}} \\ 
\left(\mathcal{M}_{\mathbf{x}_{1j}}\frac{1}{\left | \mathbf{x}-\mathbf{y} \right |} \; v_j \left(\mathbf{x}\right) \; \mathbf{e}_3 -\mathcal{M}_{\mathbf{x}_{3j}}\frac{1}{\left | \mathbf{x}-\mathbf{y} \right |} \; v_j \left(\mathbf{x}\right) \; \mathbf{e}_1\right)\cdot \mathrm{d}\boldsymbol\zeta_{\mathbf{x}} \\ 
\left(\mathcal{M}_{\mathbf{x}_{2j}}\frac{1}{\left | \mathbf{x}-\mathbf{y} \right |} \; v_j \left(\mathbf{x}\right) \; \mathbf{e}_1 -\mathcal{M}_{\mathbf{x}_{1j}}\frac{1}{\left | \mathbf{x}-\mathbf{y} \right |} \; v_j \left(\mathbf{x}\right) \; \mathbf{e}_2\right)\cdot \mathrm{d}\boldsymbol\zeta_{\mathbf{x}} 
\end{bmatrix}\mathrm{d}s_{\mathbf{y}}  \\
  &+ \frac{\mu}{4 \pi} \int\nolimits_{\Gamma}\left(\mathcal{M_\mathbf{y}\cdot\mathbf{\tilde{u}}\left ( \mathbf{y} \right ) }\right)\cdot
\int\nolimits_{\partial \Gamma}\left[\left(4\mu\mathbf{U}-\frac{2}{\left | \mathbf{x}-\mathbf{y} \right |} \mathbf{I}_3\right)\circledast \mathbf{v}\left ( \mathbf{x} \right )\right]^\top \cdot \mathrm{d}\boldsymbol\zeta_{\mathbf{x}}\;\;\mathrm{d}s_{\mathbf{y}},
\end{aligned} 
\end{equation}
where $\mathbf{I}_n$ denotes the $n\times n$ identity matrix and $\frac{\partial  }{\partial \mathbf{S} \left ( \mathbf{x} \right)}$ is the surface curl operator  $$ \frac{\partial  }{\partial \mathbf{S} \left ( \mathbf{x} \right)}  := \mathbf{n\left ( x \right )}\times \nabla_{\mathbf{x}}.$$ Here, $\otimes$ denotes the vector outer product, $\mathbf{e}_k$ is the unit vector in $k^{th}$ direction and $\epsilon_{ijk}$ is the Levi-Civita symbol. 

Similar to $\mathcal{K}$, all the differential operators in the first three surface integral terms are shifted to either $\mathbf{\tilde{u}}$ or $\mathbf{v}$. Additionally, four line integral terms arise from the right hand side of the Stokes' theorem. As with the regularized double-layer potential, the more widely recognized form of the regularized hypersingular operator (see, eg., \cite{han94}) omits these four line integrals, as they vanish for closed surfaces.

The kernels are decomposed and regularized even for the case of regular integrals, which reduces the number of required integration points in Gaussian quadrature \cite{kielhorn09}. The regularized $\mathcal{K}$ and $\mathcal{D}$ operators presented here apply to the Galerkin formulation.  However, the procedure and the formula remain identical for the collocation method, with the only difference being the outer integral. Therefore, the collocation method is not elaborated upon here. 

Although the line integrals in \eqref{newDLP} and \eqref{newHSO} involve only the fundamental solution and not its derivative, they remain strongly singular due to dimensional reduction. However, in the Galerkin formulation, no further regularization is necessary, as evidenced by the results, which will be briefly discussed in Section \ref{sec5}. 

In contrast, singularities play a role in the collocation method, and needs to be looked at closely. Considering the inner line integral of the fourth term in \eqref{newDLP}, and expanding it to the vector form yields 
\begin{equation}\label{colloexp}
\begin{aligned}
\int\nolimits_{\partial \Gamma}\frac{1}{\left | \mathbf{x}-\mathbf{y} \right |} \left ( \mathbf{\tilde{u}}\left ( \mathbf{y} \right ) \times \mathrm{d}\boldsymbol\zeta_{\mathbf{y}}\right ) = \begin{bmatrix}
 \int\nolimits_{\partial \Gamma}\frac{1}{\left | \mathbf{x}-\mathbf{y} \right |} (\tilde{u}_3\mathbf{e}_2 - \tilde{u}_2\mathbf{e}_3)\cdot \mathrm{d}\boldsymbol\zeta_{\mathbf{y}}\\
 \int\nolimits_{\partial \Gamma}\frac{1}{\left | \mathbf{x}-\mathbf{y} \right |} (\tilde{u}_1\mathbf{e}_3 - \tilde{u}_3\mathbf{e}_1)\cdot \mathrm{d}\boldsymbol\zeta_{\mathbf{y}}\\
\int\nolimits_{\partial \Gamma}\frac{1}{\left | \mathbf{x}-\mathbf{y} \right |} (\tilde{u}_2\mathbf{e}_1 - \tilde{u}_1\mathbf{e}_2)\cdot \mathrm{d}\boldsymbol\zeta_{\mathbf{y}}
\end{bmatrix}.
\end{aligned}
\end{equation}
The integrals on the right hand side are of the form
\begin{equation}\label{fpintegral}
\begin{aligned}
I(\mathbf{y})=\int\nolimits_{\partial \Gamma}\frac{1}{\left | \mathbf{x}-\mathbf{y} \right |} \;\mathbf{f(\tilde{u}(y))}\cdot \mathrm{d}\boldsymbol\zeta_{\mathbf{y}},
\end{aligned}
\end{equation}
where $I(\mathbf{y})$ does not exist as an ordinary integral when $\mathbf{x=y}$. In such cases, $I(\mathbf{y})$ must be treated as a finite-part integral \cite{paget81}, which cannot be evaluated using standard Gaussian quadrature. Integrals of the type shown in \eqref{fpintegral} are computed using the quadrature rule proposed by Paget \cite{paget81}. However, Paget's quadrature rule applies only to integrals of the form
\begin{equation}\label{Pagetintegral}
\begin{aligned}
P(f;a)=\int_{0}^{1}\frac{1}{y^a} \;f(y) \;\mathrm{d}y \;\;\;\,\;\text{for } a\geq1,
\end{aligned}
\end{equation}
and thus, the integral $I(\mathbf{y})$ must be suitably parameterized and transformed to match this structure. 

Based on the boundary triangulation \eqref{tri}, the line integral over $\partial\Gamma$ in \eqref{fpintegral} can be approximated as
\begin{equation}\label{eleLineIntegral}
\begin{aligned}
\int\nolimits_{\partial \Gamma}\frac{1}{\left | \mathbf{x}-\mathbf{y} \right |} \;\mathbf{f(\tilde{u}(y))}\cdot \mathrm{d}\boldsymbol\zeta_{\mathbf{y}} &\approx \int\nolimits_{\partial \Gamma^h}\frac{1}{\left | \mathbf{x}-\mathbf{y} \right |} \;\mathbf{f(\tilde{u}(y))}\cdot \mathrm{d}\boldsymbol\zeta_{\mathbf{y}} \\
&= \sum_{k=1}^K\sum_{i=1}^3 \int_{\gamma_i\in\tau_{k}}\frac{1}{\left | \mathbf{x}-\mathbf{y} \right |} \;\mathbf{f(\tilde{u}(y))}\cdot \mathrm{d}\boldsymbol\zeta_{\mathbf{y}},
\end{aligned}
\end{equation}
where $\gamma_i$ denotes the $i^{th}$ boundary line segment of $k^{th}$ boundary element $\tau_k$. For integration over the line segment $\gamma_i$, a parameterization with respect to $\xi$ is introduced such that,
\begin{equation}\label{LineParam}
\begin{aligned}
\mathbf{y}(\xi)= (1-\xi)\mathbf{p}_1 + \xi\mathbf{p}_2, \;\;\;\;\;\;\xi\in[0,1],
\end{aligned}
\end{equation}
where $\mathbf{p}_1$ and $\mathbf{p}_2$ are the end-points of the line segment $\gamma_i$. The differential vector term $\mathrm{d}\boldsymbol\zeta_{\mathbf{y}}$ becomes
\begin{equation}\label{diffterm}
\begin{aligned}
\mathrm{d}\boldsymbol\zeta_{\mathbf{y}} = \frac{\mathrm{d}\mathbf{y}(\xi)}{\mathrm{d}\xi}\mathrm{d}\xi = [\mathbf{p_2-p_1}]\mathrm{d}\xi. 
\end{aligned}
\end{equation}
Finally, the line integral in \eqref{eleLineIntegral} can be written as
\begin{equation}\label{paraLineIntegral}
\begin{aligned}
\int_{\gamma_i}\frac{1}{\left | \mathbf{x}-\mathbf{y} \right |} \;\mathbf{f(\tilde{u}(y))}\cdot \mathrm{d}\boldsymbol\zeta_{\mathbf{y}}=\int_0^1\frac{\mathbf{f(\xi)}}{\left | \mathbf{x}-\mathbf{y}(\xi) \right |} \;\cdot [\mathbf{p_2-p_1}]\mathrm{d}\xi,
\end{aligned}
\end{equation}
where $\mathbf{f(\xi)}:=\mathbf{f(\tilde{u}(y(\xi)))}$.

As explained in Section \ref{sec2.4}, the collocation points are chosen as the nodes of the linear triangular elements $\tau_k$. Three distinct cases of singularity are encountered, based on the position of the collocation point $\mathbf{x}$ with respect to the line segment $\gamma_i$. 
\begin{itemize}
\item \textit{Case 1: } $\mathbf{x}\notin \gamma_i$. The integral in \eqref{paraLineIntegral} is regular and is evaluated using Gaussian quadrature.
\item \textit{Case 2: } $\mathbf{x} = \mathbf{p_1}$. Substituting $\mathbf{x} = \mathbf{p_1}$ and \eqref{LineParam} into \eqref{paraLineIntegral}, the line integral becomes 
\begin{equation}\label{p1}
\begin{aligned}
\int_0^1\frac{\mathbf{f(\xi)}}{\left | \mathbf{x}-\mathbf{y}(\xi) \right |} \;\cdot [\mathbf{p_2-p_1}]\mathrm{d}\xi = \frac{1}{\left |\mathbf{p_2-p_1}\right |}\int_0^1\frac{1}{\xi} \;\mathbf{f(\xi)}\cdot [\mathbf{p_2-p_1}]\mathrm{d}\xi,
\end{aligned}
\end{equation}
which is now in the form of \eqref{Pagetintegral}. The quadrature rule from Paget \cite{paget81} can therefore be directly applied to evaluate the integral in \eqref{p1}.
\item \textit{Case 3: } $\mathbf{x} = \mathbf{p_2}$. Similar to \textit{Case 2}, the line integral becomes 
\begin{equation}\label{p2}
\begin{aligned}
\int_0^1\frac{\mathbf{f(\xi)}}{\left | \mathbf{x}-\mathbf{y}(\xi) \right |} \;\cdot [\mathbf{p_2-p_1}]\mathrm{d}\xi = \frac{1}{\left |\mathbf{p_2-p_1}\right |}\int_0^1\frac{1}{1-\xi} \;\mathbf{f(\xi)}\cdot [\mathbf{p_2-p_1}]\mathrm{d}\xi.
\end{aligned}
\end{equation}
The integral is analogous to \textit{Case 2}, but with the singularity now at $\xi=1$. With a simple change of variable, this can be transformed into the form \eqref{Pagetintegral}, allowing evaluation using Paget's quadrature rule.
\end{itemize}

All line integrals in \eqref{newDLP} are evaluated using the above approach, with the Paget's quadrature rule applied only in the collocation scheme. In the Galerkin scheme, all quadrature points lie strictly within the triangular element $\tau_k$, so only \textit{Case 1} is encountered. 
 
A direct comparison with other regularization techniques is beyond the scope of this work, but several observations can be made. Analytical integrations have the advantage of being exact, but their formulae are often complex and highly sensitive to stable implementation. Furthermore, like other expansion and cancellation methods, they are not generally applicable, as each integration depends on specific assumptions regarding geometry and shape functions. In contrast, regularization via partial integration\textemdash though requiring additional line integral evaluations in half-space problems\textemdash offers two key advantages. First, it is independent of both the problem type and the choice of shape functions. Second, nearly singular integrals become less critical, reducing the need for special handling, because the regularization applies uniformly to all integrals, including regular ones. Consequently, the number of quadrature points required for accurate evaluation is significantly reduced across all cases. 

\section{Fast multipole method}
\label{sec4}

 The BE formulations \eqref{colloMatrix} and \eqref{GalMatrix} presented in Section \ref{sec2} are solved using the biconjugate gradient stabilized method (BiCGSTAB). As stated in Section \ref{sec1}, the Chebyshev-interpolation based Fast Multipole Method (FMM) \cite{fong09} is employed to accelerate the matrix-vector multiplication and reduce both storage and computational complexity to almost linear. The black-box FMM introduced by Fong and Darve is used here with minimal modifications from \cite{fong09}.

The core idea of the FMM is to compute a low-rank approximation of the integral kernels $\mathsf{\widehat{V}}$, $\mathsf{\widehat{K}}$ and $\mathsf{\widehat{D}}$ (or $\mathsf{V}$, $\mathsf{K}$ and $\mathsf{D}$ in case of collocation) when the points $\mathbf{x}$ and $\mathbf{y}$ are sufficiently distant (the far-field). This approximation shifts the dependence of the fundamental solution $\mathbf{U}$ on the distance function $\left|\mathbf{x - y}\right|$ to two separate interpolating functions, each depending solely on either $\mathbf{x}$ or $\mathbf{y}$. For the near-field (when $\mathbf{x}$ and $\mathbf{y}$ are close), standard regularization and Gaussian quadrature are used to compute the integral kernels. To distinguish between near-field and far-field interactions, standard geometric clustering techniques from FMM \cite{liu09} are applied. This clustering is based on collocation points (or basis functions in Galerkin), but since integration occurs over the supports of the basis functions, the bounding box of each cluster is extended to fully enclose these supports. In this section, the low-rank approximation and the standard FMM operators of linear elastostatics are directly presented. A detailed description of the FMM operators, as well as the interpolation and anterpolation techniques of the multilevel scheme, can be found in \cite{greengard87} and \cite{fong09}.    

The assumption of a closed surface during the regularization process is no longer valid, as only the near-field is regularized, 
\begin{figure}[!h]
\centering
\begin{tikzpicture}
\pgfmathsetmacro{\cubex}{5}
\pgfmathsetmacro{\cubey}{2.5}
\pgfmathsetmacro{\cubez}{2.5}
\definecolor{airforceblue}{rgb}{0.36, 0.54, 0.66}

            \coordinate (A) at (0,0,0);
            \coordinate (B) at (-\cubex,0,0);
            \coordinate (C) at (-\cubex,-\cubey,0);
            \coordinate (D) at (0,-\cubey,0); 
            \coordinate (E) at (0,0,-\cubez);
            \coordinate (F) at (0,-\cubey,-\cubez);
            \coordinate (G) at (-\cubex,0,-\cubez);
            \coordinate (H) at (-\cubex*0.5,0,0);
            \coordinate (I) at (-\cubex*0.5,-\cubey,0);
            \coordinate (J) at (-\cubex*0.5,0,-\cubez);
            \coordinate (K) at (-\cubex*0.25,-\cubey*0.5,0);
            \coordinate (L) at (-\cubex*0.75,-\cubey*0.5,0);
            \coordinate (M) at (-\cubex*0.25,0,-\cubez*0.5);
            \coordinate (N) at (-\cubex*0.75,0,-\cubez*0.5);
            \coordinate (O) at (0,-\cubey*0.5,-\cubez*0.5);
            
			\draw[airforceblue, thick, fill=cyan!70] (0,0,0) -- ++(-\cubex*0.5,0,0) -- ++(0,-\cubey,0) -- ++(+\cubex*0.5,0,0) -- cycle;
			\draw[airforceblue!40, thick, fill=cyan!10] (-\cubex*0.5,0,0) -- ++(-\cubex*0.5,0,0) -- ++(0,-\cubey,0) -- ++(+\cubex*0.5,0,0) -- cycle;
			\draw[airforceblue, thick, fill=red!40] (0,0,0) -- ++(0,0,-\cubez) -- ++(0,-\cubey,0) -- ++(0,0,\cubez) -- cycle;
			\draw[airforceblue, thick, fill=cyan!70] (0,0,0) -- ++(-\cubex*0.5,0,0) -- ++(0,0,-\cubez) -- ++(\cubex*0.5,0,0) -- cycle;
			\draw[airforceblue!40, thick, fill=cyan!10] (-\cubex*0.5,0,0) -- ++(-\cubex*0.5,0,0) -- ++(0,0,-\cubez) -- ++(\cubex*0.5,0,0) -- cycle;

            
            \node at (E) [above] {{\scriptsize 1}};
            \node at (J) [above] {{\scriptsize 2}};
            \node at (G) [above] {{\scriptsize 3}};
            \node at (A) [above] {{\scriptsize 4}};
            \node at (H) [above] {{\scriptsize 5}};
            \node at (B) [left]  {{\scriptsize 6}};
            \node at (D) [below] {{\scriptsize 7}};
            \node at (I) [below] {{\scriptsize 8}};
            \node at (C) [below] {{\scriptsize 9}};
            \node at (L) [left] {{\scriptsize 10}};
            \node at (N) [left,xshift=-4pt,yshift=1pt] {{\scriptsize 11}};
            \node at (M) [left,xshift=-4pt,yshift=1pt] {{\scriptsize 12}};
            \node at (K) [left] {{\scriptsize 13}};
            \node at (F) [below] {{\scriptsize 14}};
            \node at (O) [right,xshift=-2pt] {{\scriptsize 15}};

			\draw[airforceblue, thick] (D) -- (E);
			\draw[airforceblue, thick] (A) -- (F);
			
			\draw[airforceblue, thick] (H) -- (E);
			\draw[airforceblue, thick] (J) -- (A);
			
			\draw[airforceblue!40, thick] (G) -- (H);
			\draw[airforceblue, thick] (H) -- (D);
			\draw[airforceblue!40, thick] (J) -- (B);
			
			\draw[airforceblue!40, thick] (B) -- (I);
			\draw[airforceblue, thick] (I) -- (A);
			\draw[airforceblue!40, thick] (C) -- (H);
			
			\draw[red!60, thick] (J) -- (H) -- (I);

			\fill[black!70] (A) circle (1pt);
			\fill[black!40] (B) circle (1pt);
			\fill[black!40] (C) circle (1pt);
			\fill[black!70] (D) circle (1pt);
			\fill[black!70] (E) circle (1pt);
			\fill[black!70] (F) circle (1pt);
			\fill[black!40] (G) circle (1pt);
			\fill[black!70] (H) circle (1pt);
			\fill[black!70] (I) circle (1pt);
			\fill[black!70] (J) circle (1pt);
			\fill[black!70] (K) circle (1pt);
			\fill[black!40] (L) circle (1pt);
			\fill[black!70] (M) circle (1pt);
			\fill[black!40] (N) circle (1pt);
			\fill[black!70] (O) circle (1pt);
			
\end{tikzpicture}
    \caption{Schematic representation of the near- and far-field}
    \label{fig:schematic_cuboid}
\end{figure}
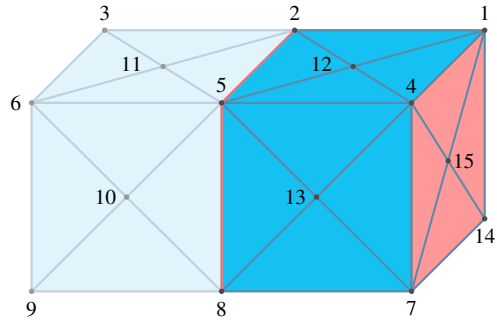
while the far-field is computed using FMM. This is illustrated schematically in Figure \ref{fig:schematic_cuboid}. Consider a cuboid geometry where the right face is subjected to Dirichlet boundary conditions, while all other edges are prescribed with Neumann boundary conditions. The figure further depicts the triangular discretization and nodal indices. To compute $\mathsf{K}$, the domain is split into a near-field (darker region) and a far-field (faded region), with node 15 being the only unknown. Since Stokes' theorem based regularization is applied exclusively to the near-field, the assumption of a closed surface is no longer valid. This necessitates computing line integrals along the boundary between the near- and far-fields (marked by red lines) to ensure valid regularization.  

An alternative approach involves applying FMM to the regularized kernel \cite{of05}, allowing for a uniform treatment of both the near- and far-field contributions while preserving the closed surface assumption. Both methods, the standard FMM  and the regularized FMM are presented here. 

\subsection{Standard FMM}
Using Chebyshev interpolation for both $\mathbf{x}$ and $\mathbf{y}$, the low-rank approximation of the fundamental solution $\mathbf{U}$ is expressed as
\begin{equation}\label{low-rankU}
\mathbf{U}\left ( \mathbf{x,y} \right ) \approx \sum_{\mathbf{n}}^{}S_{p}\left ( \mathbf{x,\bar{x}_{\mathbf{n}}} \right ) \sum_{\mathbf{m}}^{} \mathbf{U}\left ( \mathbf{\bar{x}_{\mathbf{n}},\bar{y}_{\mathbf{m}}} \right ) S_{p}\left ( \mathbf{\bar{y}_{\mathbf{m}},y} \right ).
\end{equation}
The interpolation function $S_{p}\left ( \mathbf{x,\bar{x}_{\mathbf{n}}} \right )$ is given by
\begin{equation}\label{Sp}
S_{p}\left ( \mathbf{x,\bar{x}_{\mathbf{n}}} \right ) := \prod_{i=1}^{3}S_p(x_i,\bar{x}_{n_i}),
\end{equation}
which is a three-dimensional extension of the standard one-dimensional Chebyshev interpolation function  $S_{p}\left ( x,\bar{x}_{n} \right )$, given by
\begin{equation}
S_{p}\left ( x,\bar{x}_{n} \right ) = \frac{1}{p}+\frac{2}{p}\sum_{k=1}^{p-1}T_k(x)T_k(\bar{x}_n) \;\;\;\;\;\;\;\;\;\;\; \forall \; x \in [-1,1],
\end{equation}
where $T_k(x)$ is the first-kind Chebyshev polynomial of order $k$ and $\bar{x}_{n}$ are the corresponding roots. Inserting this low-rank approximation \eqref{low-rankU} into \eqref{Vgal} gives the following expression for the discretized single layer potential $\mathsf{\widehat{V}}$,
\begin{equation}\label{VFMM}
\begin{aligned}
\mathsf{\widehat{V}}[i,j] &= \underbrace{\sum_{\mathbf{n}}\int\nolimits_{\text{supp}\left( \varphi_{i}^{0} \right)} S_{p}\left ( \mathbf{x},\bar{\mathbf{x}}_{\mathbf{n}} \right )\varphi_{i}^{0}\left(\mathbf{x}\right)}_{\textup{L2P-operator}}\: \: \underbrace{\sum_{\mathbf{m}}\mathbf{U}\left ( \bar{\mathbf{x}}_{\mathbf{n}} ,\bar{\mathbf{y}}_{\mathbf{m}} \right )}_{\textup{M2L-operator}}\: \: \\
&\cdot \underbrace{\int\nolimits_{\text{supp}\left(\varphi_{j}^{0}\right)}S_{p}\left(\bar{\mathbf{y}}_{\mathbf{m}},\mathbf{y}\right) \varphi_{j}^{0}\left(\mathbf{y}\right) }_{\textup{P2M-operator}} \; \mathrm{d} s_{\mathbf{y}}\; \mathrm{d} s_{\mathbf{x}}.
\end{aligned}
\end{equation} 
The acronyms below each term represent the standard abbreviations used in the multilevel FMM scheme. Notably, the M2L-operator, which involves the fundamental solution, needs to be evaluated only once since it depends solely on the interpolation nodes. In contrast, the other two integral kernels, involving the interpolating functions, are evaluated with the corresponding clusters of $\varphi_{i}^{0}\left(\mathbf{x}\right)$ and $\varphi_{j}^{0}\left(\mathbf{y}\right)$. Furthermore, the integral only requires the evaluation of a polynomial, rather than the entire kernel, at each quadrature point, reducing the computational effort.

For the double layer potential $\mathsf{\widehat{K}}$ and the hypersingular operator $\mathsf{\widehat{D}}$, the fundamental solution is associated with the tranction operator $\mathcal{T}$, meaning direct substitution of \eqref{low-rankU} is not feasible. However, the traction operator can be shifted onto the interpolating functions, resulting in a tensorial interpolation operator, while maintaining the same M2L-operator  \cite{schanz18}. The FMM version for $\mathsf{\widehat{K}}$ and $\mathsf{\widehat{D}}$ can then be expressed as
\begin{equation}\label{KFMM}
\begin{aligned}
\mathsf{\widehat{K}}[i,j] &= \underbrace{\sum_{\mathbf{n}}\int\nolimits_{\text{supp}\left( \varphi_{i}^{0} \right)} S_{p}\left ( \mathbf{x},\bar{\mathbf{x}}_{\mathbf{n}} \right )\varphi_{i}^{0}\left(\mathbf{x}\right)}_{\textup{L2P-operator}}\: \: \underbrace{\sum_{\mathbf{m}}\mathbf{U}\left ( \bar{\mathbf{x}}_{\mathbf{n}} ,\bar{\mathbf{y}}_{\mathbf{m}} \right )}_{\textup{M2L-operator}}\: \: \\ 
&\cdot \underbrace{\int\nolimits_{\text{supp}\left(\varphi_{j}^{1}\right)}\mathcal{T}_{\mathbf{y}}^{\top}S_{p}\left(\bar{\mathbf{y}}_{\mathbf{m}},\mathbf{y}\right) \varphi_{j}^{1}\left(\mathbf{y}\right) }_{\textup{P2M-operator}} \; \mathrm{d} s_{\mathbf{y}}\; \mathrm{d} s_{\mathbf{x}}
\end{aligned}
\end{equation} 
\begin{equation}\label{DFMM}
\begin{aligned}
\mathsf{\widehat{D}}[i,j] &= \underbrace{\sum_{\mathbf{n}}\int\nolimits_{\text{supp}\left( \varphi_{i}^{1} \right)} \mathcal{T}_{\mathbf{x}} \; S_{p}\left ( \mathbf{x},\bar{\mathbf{x}}_{\mathbf{n}} \right )\varphi_{i}^{1}\left(\mathbf{x}\right)}_{\textup{L2P-operator}}\: \: \underbrace{\sum_{\mathbf{m}}\mathbf{U}\left ( \bar{\mathbf{x}}_{\mathbf{n}} ,\bar{\mathbf{y}}_{\mathbf{m}} \right )}_{\textup{M2L-operator}}\: \: \\ 
&\cdot  \underbrace{\int\nolimits_{\text{supp}\left(\varphi_{j}^{1}\right)}\mathcal{T}_{\mathbf{y}}^{\top}S_{p}\left(\bar{\mathbf{y}}_{\mathbf{m}},\mathbf{y}\right) \varphi_{j}^{1}\left(\mathbf{y}\right) }_{\textup{P2M-operator}} \; \mathrm{d} s_{\mathbf{y}}\; \mathrm{d} s_{\mathbf{x}}
\end{aligned}
\end{equation} 
This approach ultimately requires the computation of a single M2L-operator, along with two L2P- and two P2M-operators. However, when treating the far-field using this method, it would also necessitate the computation of line integrals in \eqref{newDLP} and \eqref{newHSO} in the near field. 

\subsection{Regularized FMM}
To eliminate the need for computing line integrals, an alternative approach is to apply the low-rank approximation \eqref{low-rankU} to the regularized operators $\mathcal{K}$ \eqref{newDLP} and $\mathcal{D}$ \eqref{newHSO}, discarding the line integral terms. This allows the omission of line integrals throughout, as the entire domain can effectively be treated as a closed surface. Then, the discretized double layer potential $\mathsf{\widehat{K}}$ can be expressed as
\begin{equation}\label{regKFMM}
\begin{aligned}
\mathsf{\widehat{K}}[i,j] & = \underbrace{\sum_{\mathbf{n}}\int\nolimits_{\text{supp}\left( \varphi_{i}^{0} \right)} S_{p}\left( \mathbf{x}_{i}, \bar{\mathbf{x}}_{\mathbf{n}} \right)\varphi_{i}^{0}\left( \mathbf{x} \right)}_{\textup{L2P-operator}} \: \: \\
& \cdot \Bigg[ 2\mu \underbrace{\sum_{\mathbf{m}} \mathbf{U}\left( \bar{\mathbf{x}}_{\mathbf{n}}, \bar{\mathbf{y}}_{\mathbf{m}} \right)}_{\textup{{M2L-operator}}} \; \; \underbrace{\int\nolimits_{\text{supp}\left( \varphi_{j}^{1} \right)} S_{p}\left( \mathbf{\bar{y}}_{\mathbf{m}}, \mathbf{y} \right) \left( \mathcal{M} \varphi_{j}^{1} \right)\left( \mathbf{y} \right) \mathrm{d}s_{\mathbf{y}}}_{\textup{{P2M-SLP}}} \\
& - \frac{1}{4 \pi} \underbrace{\sum_{\mathbf{m}} \frac{1}{\left | \bar{\mathbf{x}}_{\mathbf{n}}-\bar{\mathbf{y}}_{\mathbf{m}} \right |} \; \mathbf{I}_3}_{\textup{{M2L-Laplace}}} \; \; \underbrace{\int\nolimits_{\text{supp}\left( \varphi_{j}^{1} \right)} S_{p}\left( \mathbf{\bar{y}}_{\mathbf{m}}, \mathbf{y} \right) \left( \mathcal{M} \varphi_{j}^{1} \right)\left( \mathbf{y} \right) \mathrm{d}s_{\mathbf{y}}}_{\textup{{P2M-SLP}}} \\
& + \frac{1}{4 \pi}  \underbrace{\sum_{\mathbf{m}} \frac{1}{\left | \bar{\mathbf{x}}_{\mathbf{n}}-\bar{\mathbf{y}}_{\mathbf{m}} \right |} \; \mathbf{I}_3}_{\textup{{M2L-Laplace}}} \; \; \underbrace{\int\nolimits_{\text{supp}\left( \varphi_{j}^{1} \right)} \frac{\partial}{\partial \mathbf{n\left ( y \right )}} S_{p}\left(\bar{\mathbf{y}}_{\mathbf{m}},\mathbf{y}\right) \varphi_{j}^{1}\left( \mathbf{y} \right)   \mathbf{I}_3\; \mathrm{d}s_{\mathbf{y}}}_{\textup{P2M-DLP}} \Bigg] \mathrm{d}s_{\mathbf{x}}.
\end{aligned}
\end{equation} 
The key advantage of this formulation is that it eliminates the need for computing line integrals in the near-field. However, this comes at the cost of an additional M2L-operator for the Laplacian solution, supplementing the existing one. Furthermore, two new P2M-operators, P2M-SLP and P2M-DLP, are required. The P2M-SLP applies the Günter derivative operator $\mathcal{M}$ acting on the basis functions, while the P2M-DLP operator involves applying the normal derivative to the interpolation function.   

Similarly, the discretized hypersingular operator $\mathsf{\widehat{D}}$ can be expressed as
\begin{equation}\label{regDFMM}
\begin{aligned}
\mathsf{\widehat{D}}[i,j] &=   \sum_{k=1}^{3} \Bigg[  \underbrace{\sum_{\mathbf{n}}\int\nolimits_{\text{supp}\left( \varphi_{i}^{1} \right)} S_{p}\left( \mathbf{x}_{i}, \bar{\mathbf{x}}_{\mathbf{n}} \right)\frac{\partial \varphi_{i}^{1}\left ( \mathbf{x} \right ) }{\partial S_k \left ( \mathbf{x}\right )} \; \mathbf{I}_3 }_{\textup{{L2P $\mathbb{I}_k$}}}  \\ 
 &\cdot \frac{2\mu}{4 \pi} \underbrace{\sum_{\mathbf{m}} \frac{1}{\left | \bar{\mathbf{x}}_{\mathbf{n}}-\bar{\mathbf{y}}_{\mathbf{m}} \right |} \; \mathbf{I}_3}_{\textup{{M2L-Laplace}}} \;\; \underbrace{\int\nolimits_{\text{supp}\left( \varphi_{j}^{1} \right)} S_{p}\left( \mathbf{\bar{y}}_{\mathbf{m}}, \mathbf{y} \right) \frac{\partial \varphi_{j}^{1}\left ( \mathbf{y} \right ) }{\partial S_k \left (\mathbf{y} \right )} \; \mathbf{I}_3\; \; \mathrm{d}s_{\mathbf{y}} }_{\textup{{P2M $\mathbb{I}_k$}}}\mathrm{d}s_{\mathbf{x}}\Bigg] \\ 
&- \underbrace{\sum_{\mathbf{n}}\int\nolimits_{\text{supp}\left( \varphi_{i}^{1} \right)} S_{p}\left( \mathbf{x}_{i}, \bar{\mathbf{x}}_{\mathbf{n}} \right)\mathit{diag}\left(\frac{\partial \varphi_{i}^{1}\left ( \mathbf{x} \right ) }{\partial \mathbf{S} \left ( \mathbf{x} \right )} \right )}_{\textup{{L2P $\mathbb{II}$}}}   \\ 
&\cdot \frac{\mu}{4 \pi} \underbrace{\sum_{\mathbf{m}}\frac{1}{\left | \bar{\mathbf{x}}_{\mathbf{n}}-\bar{\mathbf{y}}_{\mathbf{m}} \right |} \; \mathbf{I}_3}_{\textup{{M2L-Laplace}}} \;\; \underbrace{\int\nolimits_{\text{supp}\left( \varphi_{j}^{1} \right)} S_{p}\left( \mathbf{\bar{y}}_{\mathbf{m}}, \mathbf{y} \right) \mathbf{1}_{3} \cdot  \left[\frac{\partial \varphi_{j}^{1}\left ( \mathbf{y} \right ) }{\partial \mathbf{S} \left (\mathbf{y} \right )}\right]^{\top} \mathrm{d}s_{\mathbf{y}} }_{\textup{{P2M $\mathbb{II}$}}}\mathrm{d}s_{\mathbf{x}} \\
&-\underbrace{\sum_{\mathbf{n}}\int\nolimits_{\text{supp}\left( \varphi_{i}^{1} \right)} S_{p}\left( \mathbf{x}_{i}, \bar{\mathbf{x}}_{\mathbf{n}} \right)\left( \mathcal{M} \varphi_{i}^{1} \right)\left( \mathbf{x} \right)}_{\textup{{L2P-SLP}}}   \\
&\cdot \frac{\mu}{4 \pi} \underbrace{\sum_{\mathbf{m}} (4\mu \mathbf{U}-\frac{2}{\left | \bar{\mathbf{x}}_{\mathbf{n}}-\bar{\mathbf{y}}_{\mathbf{m}} \right |} \;\mathbf{I}_3) }_{\textup{M2L}^*} \; \; \underbrace{\int\nolimits_{\text{supp}\left( \varphi_{j}^{1} \right)} S_{p}\left( \mathbf{\bar{y}}_{\mathbf{m}}, \mathbf{y} \right) \left( \mathcal{M} \varphi_{j}^{1} \right)\left( \mathbf{y} \right) \; \; \mathrm{d}s_{\mathbf{y}} }_{\textup{{P2M-SLP}}}\mathrm{d}s_{\mathbf{x}},
\end{aligned}
\end{equation}
where $\mathbf{I}_n$ denotes the $n\times n$ identity matrix and $\mathbf{1}_{n}$ denotes a vector of ones of size $n$. As previously noted, this formulation eliminates the need for line integrals in the near-field. However, the drawback is the increased number of L2P- and P2M-operators that must be computed. Specifically, four new L2P- and P2M-operators need to be computed and stored. Nevertheless, no new M2L-operator is required, as the ones computed previously for $\mathsf{\widehat{K}}$ suffice, with M2L* being a combination of the existing operators.

Overall, the regularized FMM formulation for the discretized double layer potential $\mathsf{\widehat{K}}$ and hypersingular operator $\mathsf{\widehat{D}}$ introduces one additional M2L-operator, along with five new L2P-operators and six new P2M-operators, compared to the standard FMM. This results in more than a twofold increase in storage requirements. A comparative analysis of both FMM approaches and the role of line integrals is provided in Section \ref{sec5}.

As before, the FMM formulation presented here applies to the Galerkin BEM. Nonetheless, the methodology and expressions remain fundamentally the same for the collocation method, with the only variation arising in the integrals and the basis functions within the L2P-operators. Therefore, a detailed explanation of the collocation method is omitted.

\section{Numerical results}
\label{sec5}

This section is organized into two parts. The first part focuses on the solution of the half-space problem using both the collocation BEM and Galerkin BEM, with and without the inclusion of line integrals. The results are compared against Boussinesq's analytical solution to illustrate the impact of line integrals in half-space problems. Notably, all computations in this section are performed without employing FMM. In the second part, a convergence study is presented for two numerical examples, to evaluate and compare three different FMM formulations: the standard FMM, the standard FMM with line integrals in the near-field, and the regularized FMM. 

\subsection{Half-space problem}
The elastostatic response of a half-space subjected to a point load is considered. The half-space domain and the corresponding unbounded surface are given by 
\begin{align*}
\Omega &= \left\{ \mathbf{x} \in \mathbb{R}^3 : x_3 > 0  \right\} \\
\Gamma &= \left\{ \mathbf{x} \in \mathbb{R}^3 : x_3 = 0  \right\}.
\end{align*}
A point load $\mathbf{F}$ is applied at the origin in the positive $x_3$ direction, while the remainder of the surface is subjected to a homogeneous Neumann boundary condition. The problem statement from Section \ref{PS} is now modified to explicitly represent the Neumann problem as
\begin{align}\label{HS}
\left( \lambda + \mu \right)\nabla \left( \nabla \cdot\mathbf{u}\left ( \mathbf{x} \right ) \right) + \mu \nabla^2\mathbf{u}\left ( \mathbf{x} \right ) &= 0\; \; \;\; \; \; \; \; \;\; \; \: \forall  \:  \mathbf{x}\in\Omega \nonumber \\
\mathbf{t}\left ( \mathbf{y} \right ) &=\mathbf{F}\boldsymbol\delta(\mathbf{y})\; \; \; \; \forall  \: \mathbf{y}\in\Gamma \\
\lim_{\left| \mathbf{x} \right| \to \infty}  \left| \mathbf{x} \right|\mathbf{u(x)}&=0 \; \; \;\; \; \; \; \; \;\; \; \: \forall \; \mathbf{x}\in \overline{\Omega}. \nonumber
\end{align}
An analytical solution to this problem was first derived by Boussinesq and is commonly known as the Boussinesq-solution. This singular solution can be found in the book by Love \cite{love1892} and is given by
\begin{equation}
\begin{aligned}
u_{x_1} &= -\frac{F}{4\pi}\frac{1}{\lambda + \mu}\frac{x_1}{r}  \\
u_{x_3} &=\frac{F}{4\pi}\frac{\lambda+2\mu}{\lambda + \mu}\frac{1}{r},
\end{aligned}
\end{equation}
where $r^2 = x_1^2 + x_2^2 $ and $x_3=0$.

\begin{figure}[!h]
    \centering
    \begin{subfigure}[b]{0.45\textwidth} 
        \centering
        \includegraphics[scale=0.09]{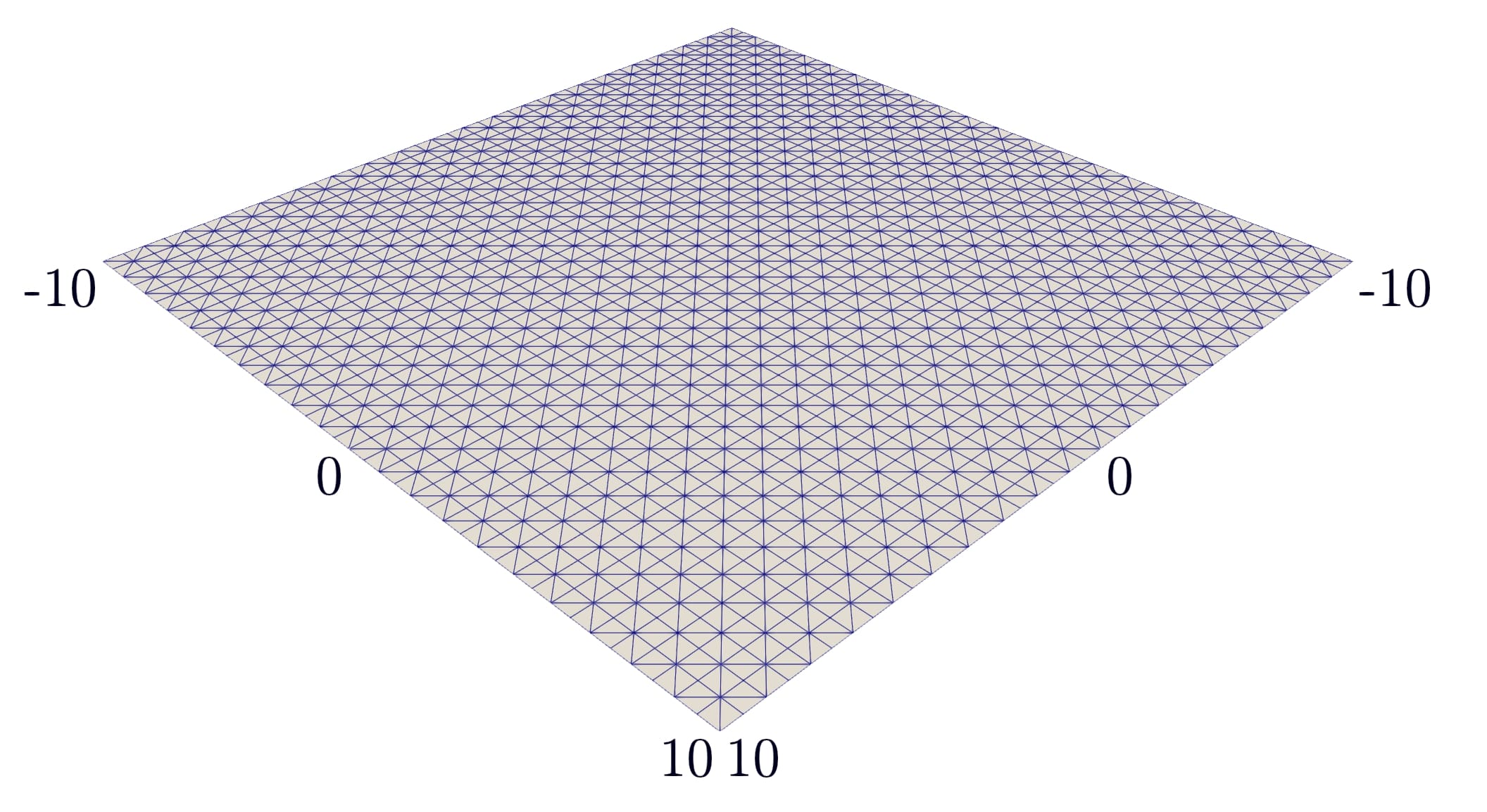}
        \caption{mesh A: 3200 elements}
        \label{fig:sheet1}
    \end{subfigure}
    \hfill
    \begin{subfigure}[b]{0.45\textwidth}
        \centering
        \includegraphics[scale=0.09]{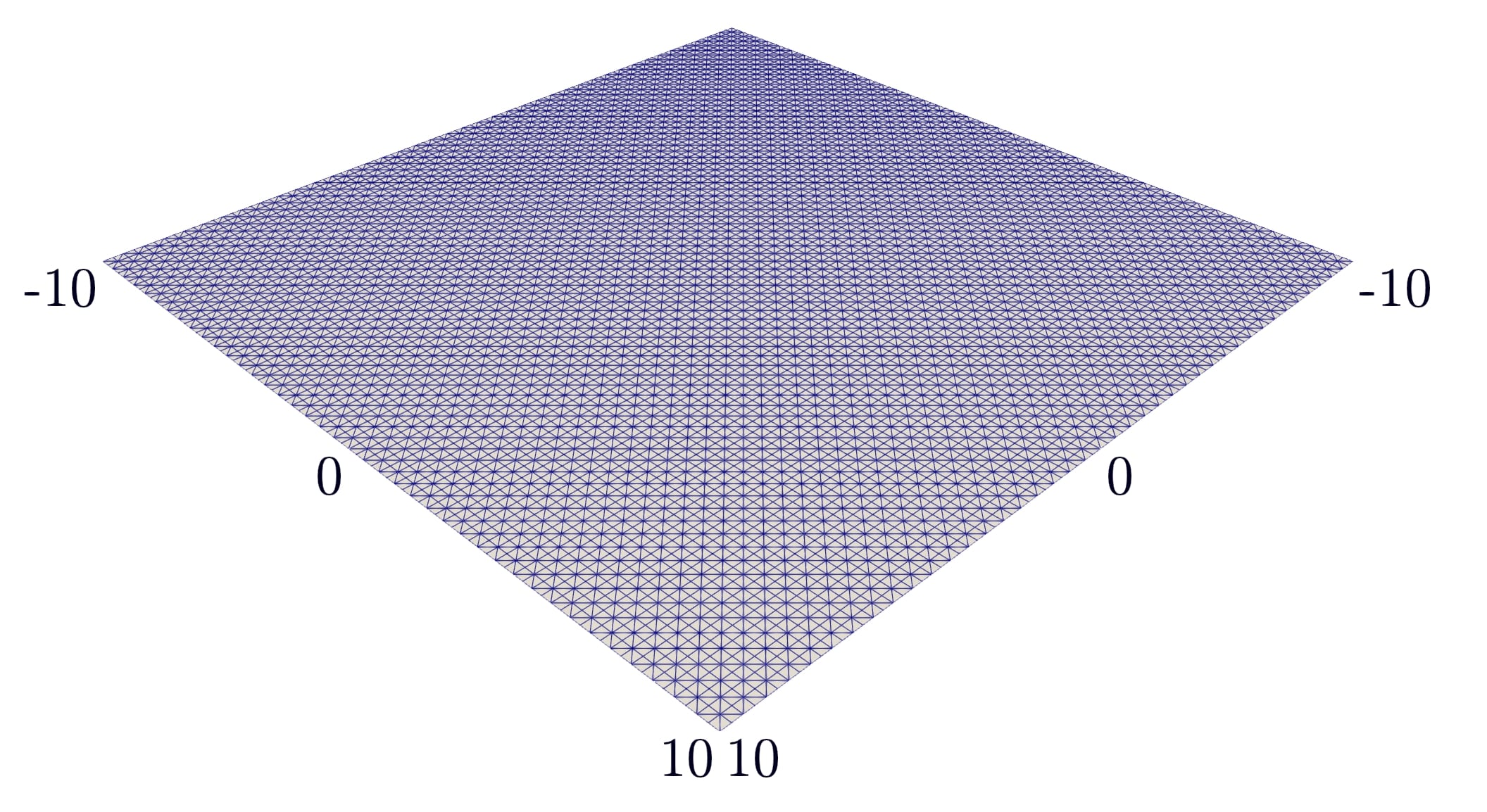}
        \caption{mesh B: 12800 elements}
        \label{fig:sheet2}
    \end{subfigure}
    \caption{Surface meshes for elastic half-space analysis}
    \label{fig:sheet}
\end{figure}

The boundary value problem formulated in \eqref{HS} is solved using both the collocation and Galerkin BEM, considering cases with and without line integral computations. The numerical results are then compared against the analytical solution. The material properties of soil are considered with $\lambda = \mu = \SI{1.3627e8}{\newton\per\meter\squared}$. To analyze the impact of discretization, two different meshes are considered, each covering a surface area of $\SI{20}{\meter}\times \SI{20}{\meter}$, as illustrated in Figure \ref{fig:sheet}.

Since a point load cannot be explicitly prescribed in the BEM framework, it is instead represented as a traction boundary condition 
\begin{equation}
\mathbf{t(y)} := 
\begin{cases}
\mathbf{F(y)} & \text{if } \mathbf{y} \in \{\Gamma : |y_1|,|y_2|\leq0.5 \; \wedge \; \|\mathbf{y}\|_1 \leq1\}, \\
\mathbf{0} & \text{everywhere else in } \Gamma
\end{cases}
\end{equation} 
on the surface $\Gamma$ where $\mathbf{F(y)}:=[0,0,1]^\top \textup{N/m}^2$. 

\begin{figure*}[!h]
    \centering
\includegraphics{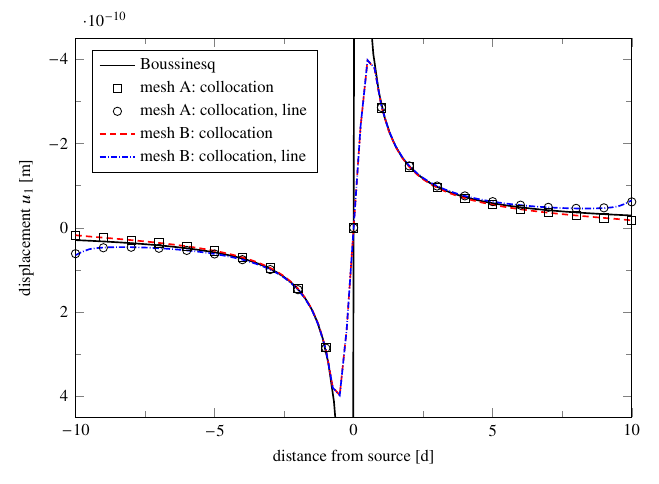}
\caption{\label{fig:Cx}Horizontal displacement with collocation}
\end{figure*}

Figure \ref{fig:Cx} shows the horizontal displacement component $u_{1}(\mathbf{y})$ along the line $L:=\{\mathbf{y}\in\Gamma:y_2=0\}$. The collocation BEM, both with and without line integrals, exhibits excellent convergence with Boussinesq’s analytical solution at most points. This aligns with the findings of \cite{schanz12}, which suggest that the collocation method performs well for reasons that remain unexplained. In this case, the inclusion of additional line integrals has little impact on the results; instead, minor deviations or kinks appear at the boundaries, likely due to singular integrals. The largest deviations between the numerical and analytical solutions occur in the region where inhomogeneous tractions are applied. This discrepancy arises because the single-point load cannot be precisely represented in the numerical model. 
\begin{figure*}[!h]
    \centering
\includegraphics{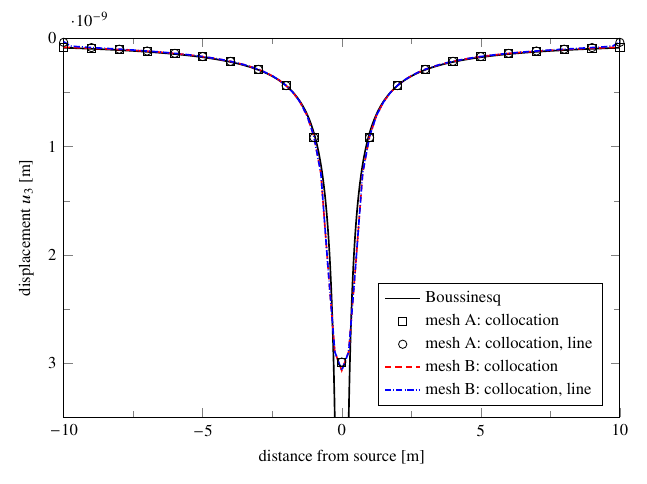}
\caption{\label{fig:Cz}Vertical displacement with collocation}
\end{figure*}
These observations also hold true for the vertical displacement $u_3(\mathbf{y})$ along the same line $L$, as shown in Figure \ref{fig:Cz}. Moreover, the results for both mesh configurations are nearly identical, exhibiting only slight variations. This further supports the conclusion that line integrals have minimal influence on the collocation approach, particularly in relation to the double layer potential.  

As previously described in Section \ref{sec3}, the collocation approach results in singular line integrals. When these integrals are evaluated using standard Gaussian quadrature, minor kinks appear near the boundary of the half-space, where the line integrals are most influential and do not cancel out.
\begin{figure*}[!h]
    \centering
\includegraphics{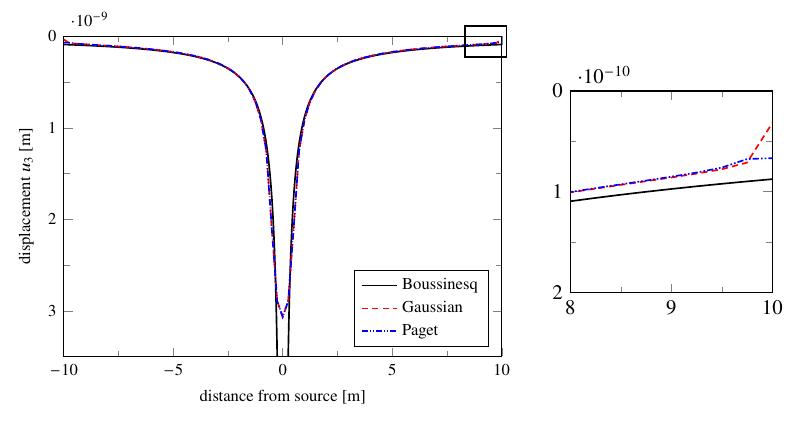}
\caption{\label{fig:Paget}Singular line integrals in collocation approach}
\end{figure*}
An example of this is clearly shown in Figure \ref{fig:Paget}, which plots the vertical displacement $u_3$, along with a zoomed-in view of the boundary region where the singular line integrals cause these kinks. However, by employing the quadrature formulas introduced by Paget \cite{paget81} without any modifications, these singularities are effectively handled, and the kinks disappear, resulting in a much smoother solution, as also demonstrated in Figure \ref{fig:Paget}. The singularity of the line integrals, combined with their minimal effect on the solution, makes their computation unnecessary in the collocation approach. 

\begin{figure*}[!h]
    \centering
\includegraphics{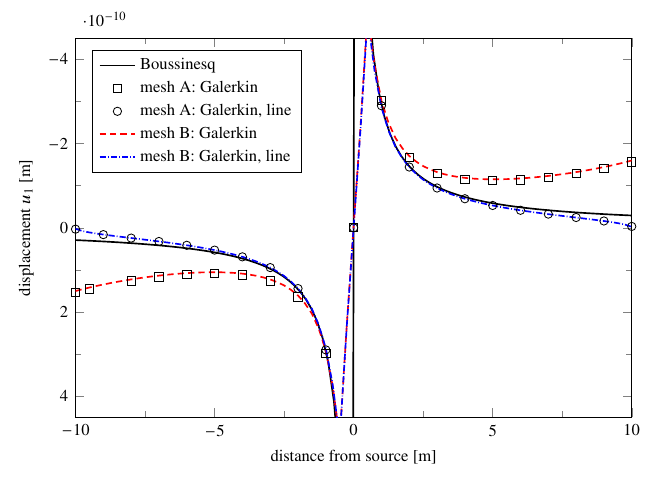}
\caption{\label{fig:Gx}Horizontal displacement with Galerkin}
\end{figure*}

\begin{figure*}[!h]
    \centering
\includegraphics{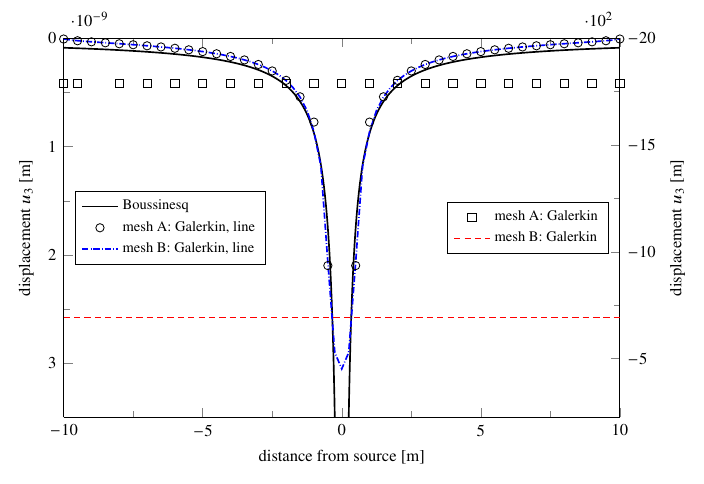}
\caption{\label{fig:Gz}Vertical displacement with Galerkin}
\end{figure*}

The Galerkin approach presents a different scenario altogether. Figure \ref{fig:Gx} shows the horizontal displacement component $u_{1}(\mathbf{y})$ along the line $L := \{\mathbf{y} \in \Gamma : y_2 = 0\}$ for the Galerkin approach.
In this case, the numerical results obtained from the Galerkin BEM without the additional line integrals show significant deviation from the analytical solution for both discretization meshes A and B. However, when line integrals are introduced into the hypersingular operator \eqref{newHSO}, the results improve substantially, converging very closely to Boussinesq's analytical solution. This demonstrates the importance of line integrals, which is further emphasized in the vertical displacement $u_{3}(\mathbf{y})$, shown in Figure \ref{fig:Gz}.

To better illustrate the extent of these differences, two y-axes with different scales are presented. The y-axis on the right corresponds to the Galerkin approach without line integrals, while the left-hand side shows the analytical solution alongside the Galerkin approach with line integrals.  Without the inclusion of the additional line integrals, the Galerkin approach fails completely, producing garbage values. However, upon incorporating the four line integral terms from \eqref{newHSO}, the solution significantly improves and converges closely to the analytical solution, except in the region near the inhomogeneous traction boundary, where the discrepancy is expected.

These results clearly demonstrate that, in the context of half-space problems, the collocation method performs surprisingly well without the need for line integrals. On the other hand, the Galerkin approach involving the hypersingular operator, is much more sensitive and fails completely without the line integral computations. Furthermore, the line integrals in the Galerkin approach do not exhibit any singularities, as described in Section \ref{sec3} and evidenced by the smoothness of the solution, which lacks any irregularities or kinks. 

\subsection{FMM}
A convergence study is conducted on two numerical examples, both prescribed with mixed boundary conditions \eqref{BCs}. The first example is solved using the collocation approach, while the second employs the Galerkin approach. In both cases, the material parameters are set as $\lambda = \SI{0.2778}{\newton\per\meter\squared}$ and $\mu = \SI{0.4167}{\newton\per\meter\squared}$. The boundary conditions are imposed using the fundamental solution as 
\begin{equation}
\begin{aligned}
\mathbf{g}_{D}\left ( \mathbf{x} \right )&:=\hat{\mathbf{u}}\left ( \mathbf{x} \right )=\mathbf{U} \left( \mathbf{x,y_\mathit{S}}\right) \; \; \; \; \; \; \; \; \;\; \;\; \; \:    \mathbf{x}\in\Gamma \; \; \; \; \mathbf{y_\mathit{S}}\notin\Omega \\
\mathbf{g}_{N}\left ( \mathbf{x} \right )&:=\hat{\mathbf{t}}\left ( \mathbf{x} \right )=\left ( \mathcal{T}_{\mathbf{y}}\mathbf{U} \right )^{\top}\left( \mathbf{x,y_\mathit{S}}\right) \; \; \; \; \; \; \mathbf{x}\in\Gamma,
\end{aligned}
\end{equation}
where $\mathbf{y_\mathit{S}}$ is the source point located outside the computational domain.  Since $\mathbf{U}$ is an analytic solution of the elastostatic equation, the relative $L_2$ errors are computed as
\begin{equation}
\begin{aligned}
\textup{err}_{\mathit{rel}}\left ( \mathbf{u} \right )&=\left\|\hat{\mathbf{u}}-\tilde{\mathbf{u}} \right\|_{L_{2}\left ( \Gamma  \right )}\\
\textup{err}_{\mathit{rel}}\left ( \mathbf{t} \right )&=\left\|\hat{\mathbf{t}}-\tilde{\mathbf{t}} \right\|_{L_{2}\left ( \Gamma  \right )}. 
\end{aligned}
\end{equation}
The order of convergence (eoc) is given by 
\begin{equation}
\begin{aligned}
\text{eoc} = \text{log}_2\left(\frac{\text{err}_{rel}^k}{\text{err}_{rel}^{k+1}}\right),
\end{aligned}
\end{equation}
where $k$ denotes the corresponding refinement level. Based on the spatial discretization, it is expected that the Neumann data exhibits a linear order of convergence, while the Dirichlet data exhibits a quadratic order of error convergence \cite{steinbach08}.

\begin{figure}[!h]
    \centering
    \begin{subfigure}[b]{0.45\textwidth} 
        \centering
        \includegraphics[scale=0.1]{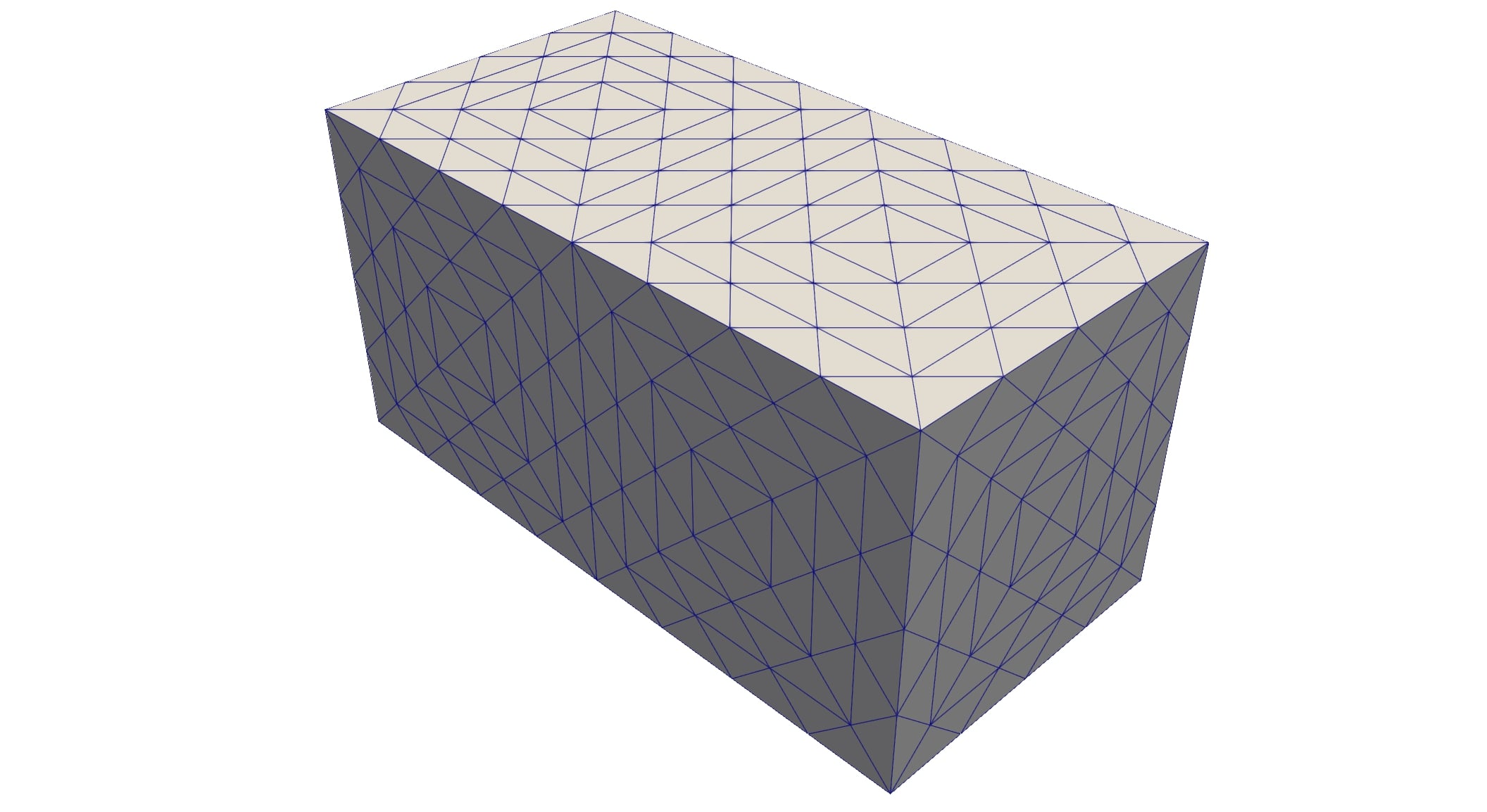}
        \caption{lvl 2: 640 elements}
        \label{fig:cuboid1}
    \end{subfigure}
    \begin{subfigure}[b]{0.45\textwidth}
        \centering
        \includegraphics[scale=0.1]{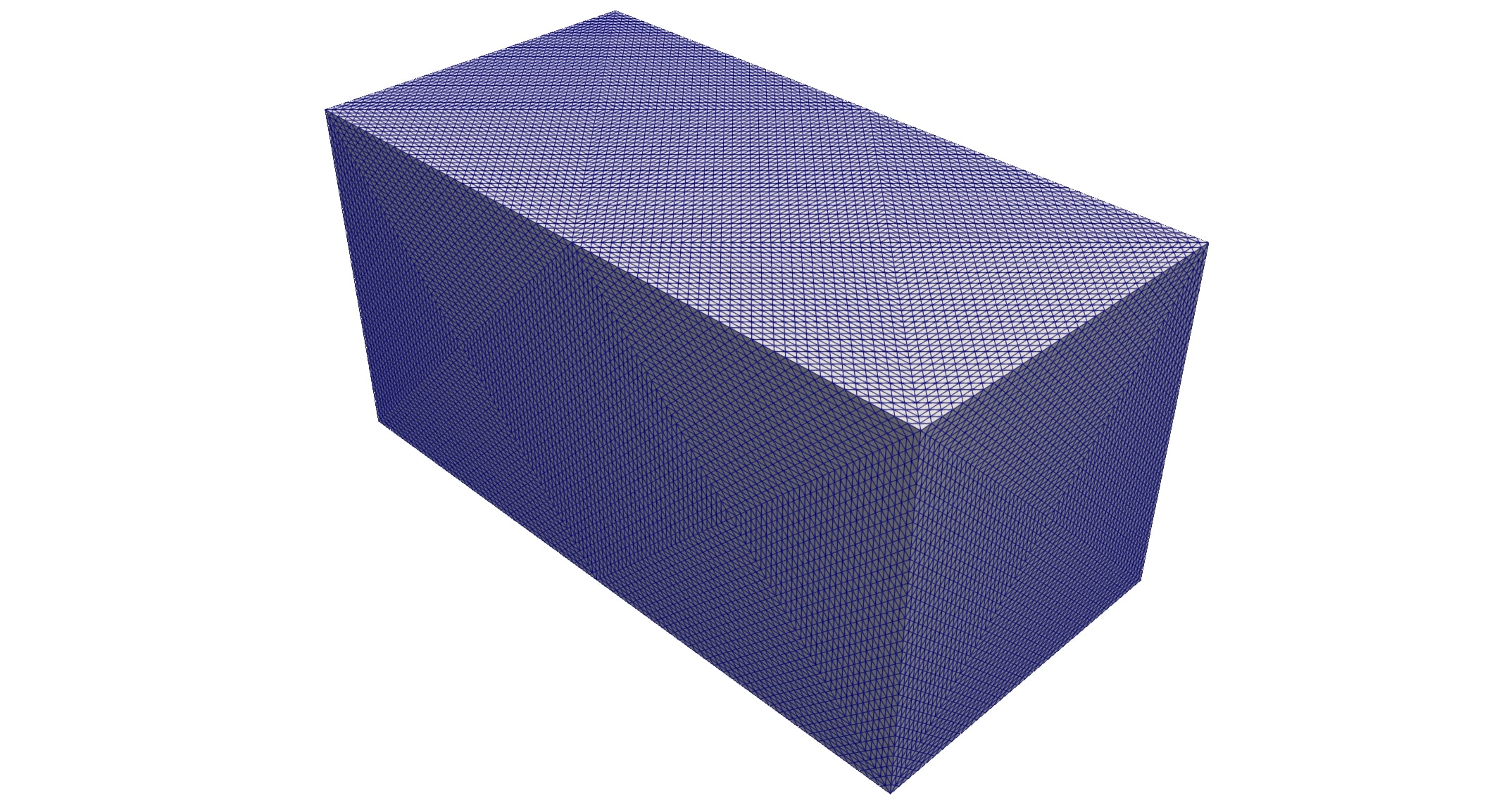}
        \caption{lvl 5: 40960 elements}
        \label{fig:cuboid2}
    \end{subfigure}
    \caption{Cuboid}
    \label{fig:cuboid}
\end{figure}

The first example considers a cuboid of dimensions $\SI{2}{\meter} \times \SI{1}{\meter} \times \SI{1}{\meter}$, as shown in Figure \ref{fig:cuboid}, with the coordinate system located at the bottom-left corner of its right face.
\begin{table}[!ht]
    \centering
    \renewcommand{\arraystretch}{1.25}
    \begin{tabular}{crrrr}
    \hline\hline
       lvl & dof & $h$ & $p$ & ${\text{F}}_{lvl}$  \\
        \hline
       0 & 63 & $\SI{1.0}{\meter}$ & 2 & 1 \\
       1 & 255 & $\SI{0.5}{\meter}$ & 3 & 2 \\
       2 & 1035 & $\SI{0.25}{\meter}$ & 4 & 2 \\
       3 & 4179 & $\SI{0.125}{\meter}$ & 5 & 3 \\
       4 & 16803 & $\SI{0.0625}{\meter}$ & 6 & 4 \\
       5 & 67395 & $\SI{0.03125}{\meter}$ & 7 & 5 \\
        \hline\hline
    \end{tabular}
    \caption{Discretization and FMM parameters}
    \label{table: ref levels 1}
\end{table}
The initial unrefined mesh (lvl 0) consists of 40 elements, with 5 levels of refinement, and the finest mesh (lvl 5) contains 40960 elements. Dirichlet boundary conditions are applied to the right face of the cuboid, while Neumann boundary conditions are imposed on the remaining faces. The source point is  located at $\mathbf{y_\mathit{S}} = (3.0,0.5,0.5)^\top$. Table \ref{table: ref levels 1} presents the discretization parameters along with the corresponding order of interpolation ($p$) and levels (${\text{F}}_{lvl}$) of the multi-level FMM scheme, at each refinement level.

This problem is solved using the collocation approach with three variations of FMM: the standard FMM (FMM), standard FMM with line integrals in the
\begin{table}[!ht]
    \centering
    \begin{minipage}{0.45\textwidth} 
        \centering
        \renewcommand{\arraystretch}{1}
        \begin{tabular}{clll}
            \hline\hline
            lvl & FMM & FMML & RFMM \\
            \hline
            0 & 6.80e-1 & 6.80e-1 & 6.80e-1 \\
			  & X       & X       & X    	\\     
			  \hline
			  \addlinespace[0.5ex]      
            1 & 2.15e-1 & 2.15e-1 & 2.15e-1 \\
              & 1.661   & 1.661   & 1.661   \\ 
              \hline
              \addlinespace[0.5ex]
            2 & 5.07e-2 & 5.07e-2 & 5.07e-2 \\
              & 2.084   & 2.084   & 2.085   \\ 
              \hline
              \addlinespace[0.5ex]
            3 & 1.10e-2 & 1.10e-2 & 1.16e-2 \\
              & 2.208   & 2.208   & 2.118   \\ 
              \hline
              \addlinespace[0.5ex]
            4 & 2.36e-3 & 2.36e-3 & 2.57e-3 \\
              & 2.219   & 2.219   & 2.180   \\ 
              \hline
              \addlinespace[0.5ex]
            5 & 6.62e-4 & 6.62e-4 & 9.26e-4 \\
              & 1.832   & 1.832   & 1.477   \\ 
            \hline\hline
        \end{tabular}
        \caption{$\textup{err}_{\mathit{rel}}\left ( \mathbf{u} \right )$ and eoc}
        \label{table:error1}
    \end{minipage}%
    \hspace{0.5cm} 
    \begin{minipage}{0.45\textwidth} 
        \centering
        \renewcommand{\arraystretch}{1}
        \begin{tabular}{clll}
            \hline\hline
            lvl & FMM & FMML & RFMM \\
            \hline
            0 & 4.52e-1 & 4.52e-1 & 4.52e-1 \\
			  & X       & X       & X    	\\    
			   \hline
			  \addlinespace[0.5ex] 
            1 & 2.08e-1 & 2.08e-1 & 2.08e-1 \\
              & 1.119   & 1.119   & 1.119   \\ 
              \hline
              \addlinespace[0.5ex]
            2 & 9.37e-2 & 9.37e-2 & 9.37e-2 \\
              & 1.153   & 1.153   & 1.153   \\ 
              \hline
              \addlinespace[0.5ex]
            3 & 4.55e-2 & 4.55e-2 & 4.55e-2 \\
              & 1.041   & 1.041   & 1.042   \\ 
              \hline
              \addlinespace[0.5ex]
            4 & 2.26e-2 & 2.26e-2 & 2.26e-2 \\
              & 1.011   & 1.011   & 1.010   \\ 
              \hline
              \addlinespace[0.5ex]
            5 & 1.13e-2 & 1.13e-2 & 1.13e-2 \\
              & 1.001   & 1.001   & 1.001   \\ 
            \hline\hline
        \end{tabular}
        \caption{$\textup{err}_{\mathit{rel}}\left ( \mathbf{t} \right )$ and eoc}
        \label{table:error2}
    \end{minipage}
\end{table}
near-field (FMML), and regularized FMM (RFMM). Relative errors (above) and the estimated order of convergence (eoc) (below) are summarized in Tables \ref{table:error1} and \ref{table:error2}. The Dirichlet errors for FMM and FMML are identical, suggesting that the line integrals do not have any effect on the solution. The errors for RFMM are also similar, with only minor variations at finer refinerment levels, indicating that the regularization of the far-field does not enhance or modify the solution. For the Neumann errors, all three versions of FMM yield identical errors as expected, since the single-layer potential remains unchanged. For both solutions, the eoc falls within the expected range.

The second example involves a cube with dimensions  $\SI{1}{\meter} \times \SI{1}{\meter} \times \SI{1}{\meter}$, from which a smaller cube of size ($\SI{0.5}{\meter} \times \SI{0.5}{\meter} \times \SI{0.5}{\meter}$) is removed. The removed portion is located at one of the cube's corners, effectively creating an L-shaped geometry with a void at the corner, often referred to as the Fichera cube. The coordinate system is located at the center of the cube. The problem is solved using the Galerkin approach.
\begin{figure}[!h]
    \centering
    \begin{subfigure}[b]{0.45\textwidth} 
        \centering
        \includegraphics[scale=0.1]{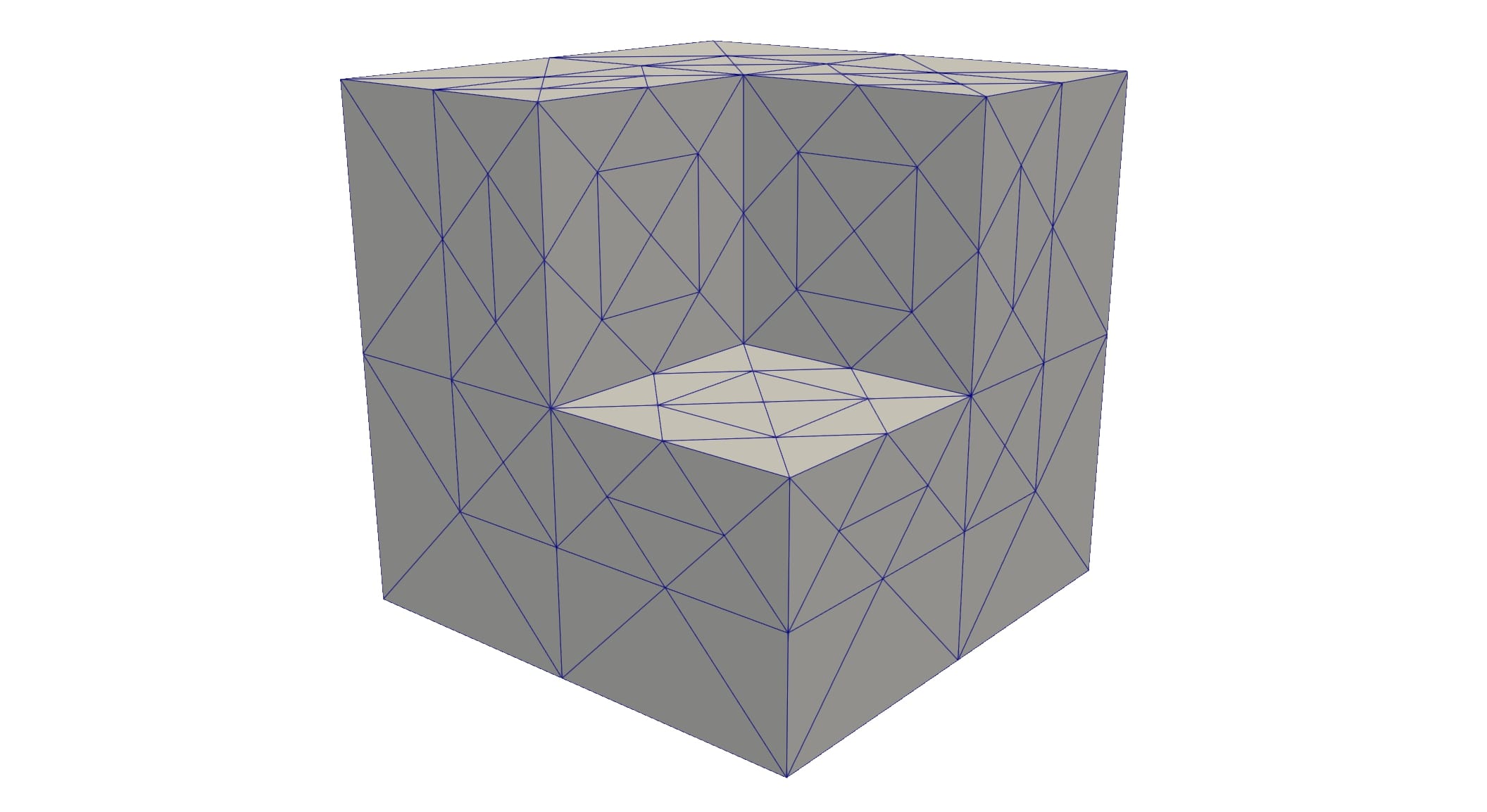}
        \caption{lvl 0: 186 elements}
        \label{fig:lcube1}
    \end{subfigure}
    \begin{subfigure}[b]{0.45\textwidth}
        \centering
        \includegraphics[scale=0.1]{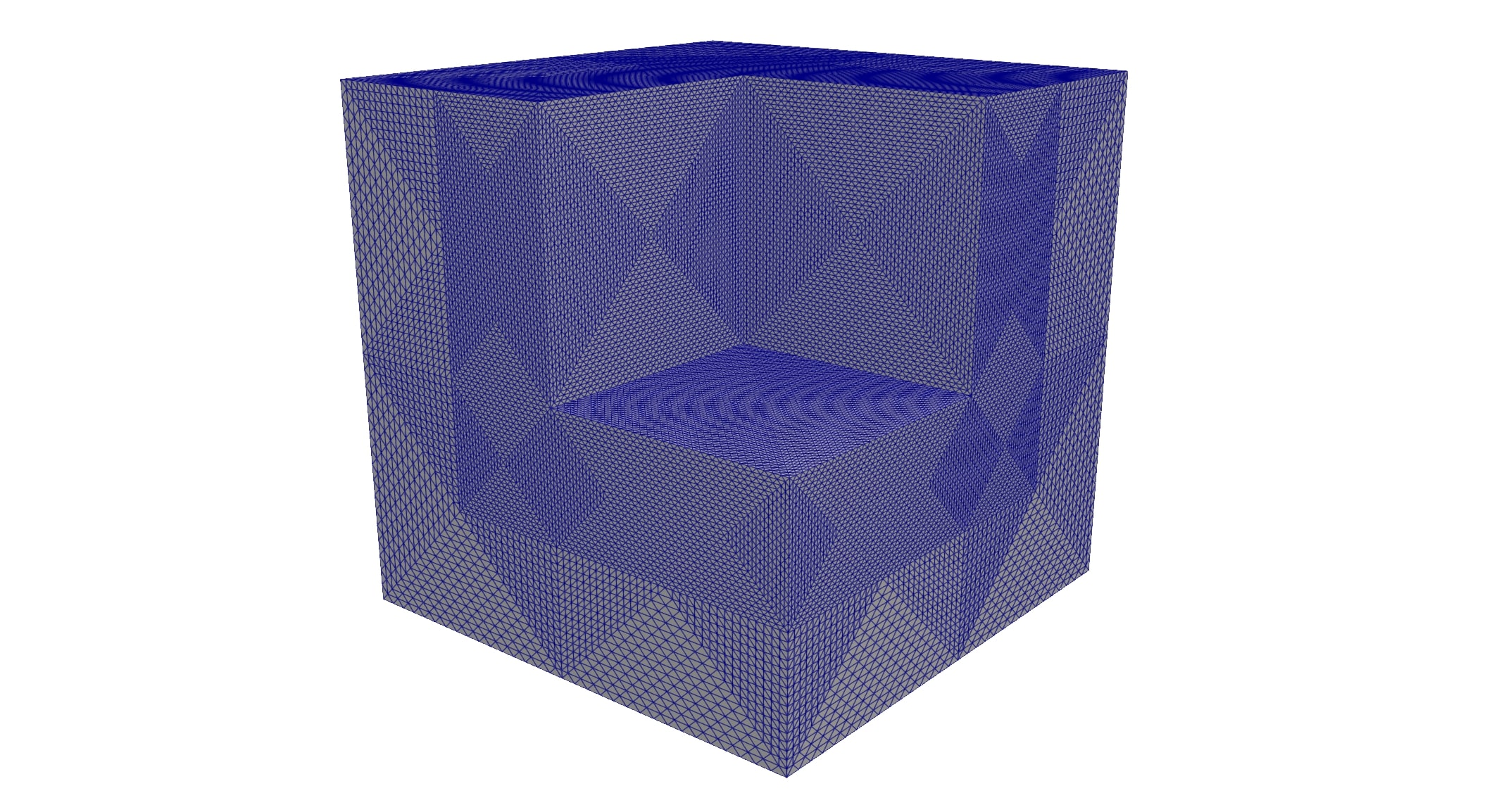}
        \caption{lvl 4: 47616 elements}
        \label{fig:lcube2}
    \end{subfigure}
    \caption{Fichera cube}
    \label{fig:lcube}
\end{figure}

The initial unrefined mesh (lvl 0) consists of 186 elements, with four levels of refinement, 
\begin{table}[!h]
    \centering
    \renewcommand{\arraystretch}{1.25}
    \begin{tabular}{crrrr}
    \hline\hline
       lvl & dof & $h$ & $p$ & ${\text{F}}_{lvl}$  \\
        \hline
       0 & 498 & $\SI{0.5}{\meter}$ & 2 & 1 \\
       1 & 2019 & $\SI{0.25}{\meter}$ & 3 & 2 \\
       2 & 8139 & $\SI{0.125}{\meter}$ & 4 & 2 \\
       3 & 32691 & $\SI{0.0625}{\meter}$ & 5 & 3 \\
       4 & 131043 & $\SI{0.03125}{\meter}$ & 6 & 4 \\
        \hline\hline
    \end{tabular}
    \caption{Discretization and FMM parameters}
    \label{table: ref levels 2}
\end{table}
leading to a finest mesh (lvl 4) containing 47616 elements, as shown in Figure \ref{fig:lcube}. Neumann boundary conditions are applied to the front-left face of the Fichera cube, while Dirichlet boundary conditions are prescribed on the remaining faces. The source point is  located at $\mathbf{y_\mathit{S}} = (0.75,-0.15,0.1)^\top$. The discretization parameters along with the interpolation order ($p$) and levels (${\text{F}}_{lvl}$) of the multi-level FMM scheme for different refinement levels are listed in Table \ref{table: ref levels 2}. 

Once again, three versions of FMM are used for the convergence study, and the corresponding  relative $L_2$ errors (above) and eoc (below) are presented in Tables \ref{table:error3} and \ref{table:error4}.
\begin{table}[!h]
    \centering
    \begin{minipage}{0.45\textwidth} 
        \centering
        \renewcommand{\arraystretch}{1.25}
        \begin{tabular}{clll}
            \hline\hline
            lvl & FMM & FMML & RFMM \\
            \hline
            0 & 4.18e-2 & 4.18e-2 & 4.18e-2 \\
			  & X       & X       & X    	\\    
			   \hline
			  \addlinespace[0.5ex] 
            1 & 1.17e-2 & 1.17e-2 & 1.17e-2 \\
              & 1.836   & 1.836   & 1.836   \\ 
              \hline
              \addlinespace[0.5ex]
            2 & 3.06e-3 & 3.06e-3 & 3.06e-3 \\
              & 1.934   & 1.934   & 1.935   \\ 
              \hline
              \addlinespace[0.5ex]
            3 & 7.83e-4 & 7.83e-4 & 7.75e-4 \\
              & 1.968   & 1.968   & 1.982   \\ 
              \hline
              \addlinespace[0.5ex]
            4 & 2.14e-4 & 2.14e-4 & 2.04e-4 \\
              & 1.869   & 1.869   & 1.929   \\ 
            \hline\hline
        \end{tabular}
        \caption{$\textup{err}_{\mathit{rel}}\left ( \mathbf{u} \right )$ and eoc}
        \label{table:error3}
    \end{minipage}%
    \hspace{0.5cm} 
    \begin{minipage}{0.45\textwidth} 
        \centering
        \renewcommand{\arraystretch}{1.25}
        \begin{tabular}{clll}
            \hline\hline
            lvl & FMM & FMML & RFMM \\
            \hline
            0 & 4.66e-1 & 4.66e-1 & 4.66e-1 \\
			  & X       & X       & X    	\\    
			   \hline
			  \addlinespace[0.5ex] 
            1 & 2.27e-1 & 2.27e-1 & 2.27e-1 \\
              & 1.037   & 1.037   & 1.037   \\ 
              \hline
              \addlinespace[0.5ex]
            2 & 1.07e-1 & 1.07e-1 & 1.07e-1 \\
              & 1.091   & 1.091   & 1.091   \\ 
              \hline
              \addlinespace[0.5ex]
            3 & 5.09e-2 & 5.09e-2 & 5.08e-2 \\
              & 1.066   & 1.066   & 1.070   \\ 
              \hline
              \addlinespace[0.5ex]
            4 & 2.52e-2 & 2.52e-2 & 2.51e-2 \\
              & 1.014   & 1.014   & 1.018   \\ 
            \hline\hline
        \end{tabular}
        \caption{$\textup{err}_{\mathit{rel}}\left ( \mathbf{t} \right )$ and eoc}
        \label{table:error4}
    \end{minipage}
\end{table}
The results align with those of the collocation approach, where the identical Dirichlet errors for FMM and FMML confirm that line integrals have no impact on the solution. This implies that, despite the assumption of a closed surface being invalid, the line integrals effectively cancel out. The RFMM errors follow the same trend, differing only slightly at finer refinement levels, suggesting that far-field regularization neither improves nor alters the solution. As before, the Neumann errors remain unchanged for all FMM versions and the eoc for both the errors falls within the expected range.

The main drawback of RFMM, despite achieving similar results, is its higher storage requirement. This is due to the additional operators that need to be computed. 
\begin{table}[!h]
    \centering
        \renewcommand{\arraystretch}{1.25}
    \begin{tabular}{c|ll|ll}
        \hline\hline
        lvl & \multicolumn{2}{c|}{storage [\SI{}{\giga\byte}]} & \multicolumn{2}{c}{time [\SI{}{\second}]} \\
        
             & FMM & RFMM & FMM & RFMM \\
        \hline
        2    & 3.48e-3    & 6.65e-3     & 0.2   & 0.3    \\
        3    & 2.08e-1    & 4.17e-1     & 15.8   & 19.7    \\
        4    & 8.13e-1    & 1.61     & 338.6   & 344.5    \\
        5    & 2.18    & 4.23     & 2428.1   &  2745.6   \\
        \hline\hline
    \end{tabular}
        \caption{Computational costs}
        \label{table:comp1}
\end{table}

Table \ref{table:comp1} shows the storage requirements and the pre-computing time for both the standard FMM and RFMM, for the first example problem. The first two refinement levels are omitted because FMM is not triggered for those levels, as the admissibility condition for the far-field is not satisfied. However, for the finer refinement levels, the storage requirements for the RFMM is approximately doubled, mainly due to the fact that an additional M2L operator is computed. The pre-computation time for RFMM is also slightly higher, indicating that RFMM incurs extra computational costs to achieve similar results.

For the FMML, the cancellation of the line integrals is illustrated with the same schematic representation as in Figure \ref{fig:schematic_cuboid}.
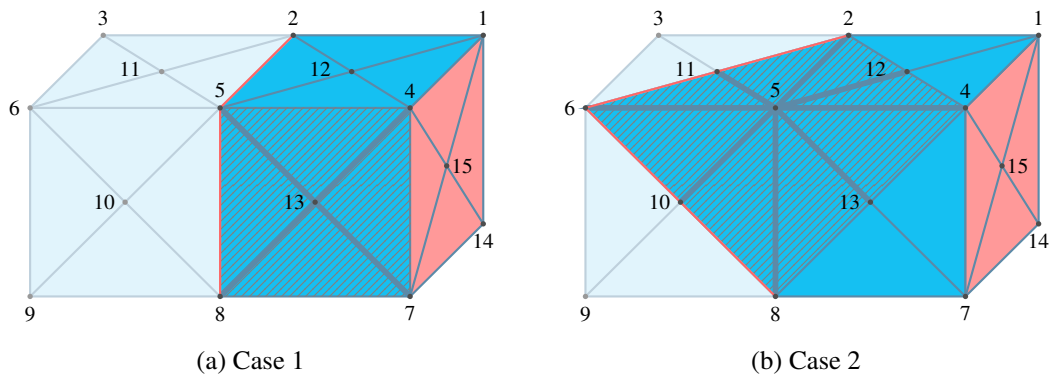
\begin{figure}[!h]

    \begin{subfigure}[b]{0.45\textwidth} 

\begin{tikzpicture}
\pgfmathsetmacro{\cubex}{5}
\pgfmathsetmacro{\cubey}{2.5}
\pgfmathsetmacro{\cubez}{2.5}
\definecolor{airforceblue}{rgb}{0.36, 0.54, 0.66}

            \coordinate (A) at (0,0,0);d
            \coordinate (B) at (-\cubex,0,0);
            \coordinate (C) at (-\cubex,-\cubey,0);
            \coordinate (D) at (0,-\cubey,0); 
            \coordinate (E) at (0,0,-\cubez);
            \coordinate (F) at (0,-\cubey,-\cubez);
            \coordinate (G) at (-\cubex,0,-\cubez);
            \coordinate (H) at (-\cubex*0.5,0,0);
            \coordinate (I) at (-\cubex*0.5,-\cubey,0);
            \coordinate (J) at (-\cubex*0.5,0,-\cubez);
            \coordinate (K) at (-\cubex*0.25,-\cubey*0.5,0);
            \coordinate (L) at (-\cubex*0.75,-\cubey*0.5,0);
            \coordinate (M) at (-\cubex*0.25,0,-\cubez*0.5);
            \coordinate (N) at (-\cubex*0.75,0,-\cubez*0.5);
            \coordinate (O) at (0,-\cubey*0.5,-\cubez*0.5);
            
			\draw[airforceblue, thick, fill=cyan!70] (0,0,0) -- ++(-\cubex*0.5,0,0) -- ++(0,-\cubey,0) -- ++(+\cubex*0.5,0,0) -- cycle;
			\draw[airforceblue!40, thick, fill=cyan!10] (-\cubex*0.5,0,0) -- ++(-\cubex*0.5,0,0) -- ++(0,-\cubey,0) -- ++(+\cubex*0.5,0,0) -- cycle;
			\draw[airforceblue, thick, fill=red!40] (0,0,0) -- ++(0,0,-\cubez) -- ++(0,-\cubey,0) -- ++(0,0,\cubez) -- cycle;
			\draw[airforceblue, thick, fill=cyan!70] (0,0,0) -- ++(-\cubex*0.5,0,0) -- ++(0,0,-\cubez) -- ++(\cubex*0.5,0,0) -- cycle;
			\draw[airforceblue!40, thick, fill=cyan!10] (-\cubex*0.5,0,0) -- ++(-\cubex*0.5,0,0) -- ++(0,0,-\cubez) -- ++(\cubex*0.5,0,0) -- cycle;
			
			\usetikzlibrary{patterns}
			\fill[pattern=north east lines, thick, pattern color=gray] (A) -- (H) -- (I) -- (D) -- cycle;
                 
            
            \node at (E) [above] {{\scriptsize 1}};
            \node at (J) [above] {{\scriptsize 2}};
            \node at (G) [above] {{\scriptsize 3}};
            \node at (A) [above] {{\scriptsize 4}};
            \node at (H) [above] {{\scriptsize 5}};
            \node at (B) [left]  {{\scriptsize 6}};
            \node at (D) [below] {{\scriptsize 7}};
            \node at (I) [below] {{\scriptsize 8}};
            \node at (C) [below] {{\scriptsize 9}};
            \node at (L) [left] {{\scriptsize 10}};
            \node at (N) [left,xshift=-4pt,yshift=1pt] {{\scriptsize 11}};
            \node at (M) [left,xshift=-4pt,yshift=1pt] {{\scriptsize 12}};
            \node at (K) [left] {{\scriptsize 13}};
            \node at (F) [below] {{\scriptsize 14}};
            \node at (O) [right,xshift=-2pt] {{\scriptsize 15}};

			\draw[airforceblue, thick] (D) -- (E);
			\draw[airforceblue, thick] (A) -- (F);
			
			\draw[airforceblue, thick] (H) -- (E);
			\draw[airforceblue, thick] (J) -- (A);
			
			\draw[airforceblue!40, thick] (G) -- (H);
			\draw[airforceblue, line width=0.75mm] (H) -- (D);
			\draw[airforceblue!40, thick] (J) -- (B);
			
			\draw[airforceblue!40, thick] (B) -- (I);
			\draw[airforceblue, line width=0.75mm] (I) -- (A);
			\draw[airforceblue!40, thick] (C) -- (H);
			
			\draw[red!60, thick] (J) -- (H) -- (I);

			\fill[black!70] (A) circle (1pt);
			\fill[black!40] (B) circle (1pt);
			\fill[black!40] (C) circle (1pt);
			\fill[black!70] (D) circle (1pt);
			\fill[black!70] (E) circle (1pt);
			\fill[black!70] (F) circle (1pt);
			\fill[black!40] (G) circle (1pt);
			\fill[black!70] (H) circle (1pt);
			\fill[black!70] (I) circle (1pt);
			\fill[black!70] (J) circle (1pt);
			\fill[black!70] (K) circle (1pt);
			\fill[black!40] (L) circle (1pt);
			\fill[black!70] (M) circle (1pt);
			\fill[black!40] (N) circle (1pt);
			\fill[black!70] (O) circle (1pt);
			
\end{tikzpicture}
    \caption{Case 1}
    \label{fig:case1}
    \end{subfigure}
    \hspace{5mm}
    \begin{subfigure}[b]{0.45\textwidth}
\begin{tikzpicture}
\pgfmathsetmacro{\cubex}{5}
\pgfmathsetmacro{\cubey}{2.5}
\pgfmathsetmacro{\cubez}{2.5}
\definecolor{airforceblue}{rgb}{0.36, 0.54, 0.66}

            \coordinate (A) at (0,0,0);
            \coordinate (B) at (-\cubex,0,0);
            \coordinate (C) at (-\cubex,-\cubey,0);
            \coordinate (D) at (0,-\cubey,0); 
            \coordinate (E) at (0,0,-\cubez);
            \coordinate (F) at (0,-\cubey,-\cubez);
            \coordinate (G) at (-\cubex,0,-\cubez);
            \coordinate (H) at (-\cubex*0.5,0,0);
            \coordinate (I) at (-\cubex*0.5,-\cubey,0);
            \coordinate (J) at (-\cubex*0.5,0,-\cubez);
            \coordinate (K) at (-\cubex*0.25,-\cubey*0.5,0);
            \coordinate (L) at (-\cubex*0.75,-\cubey*0.5,0);
            \coordinate (M) at (-\cubex*0.25,0,-\cubez*0.5);
            \coordinate (N) at (-\cubex*0.75,0,-\cubez*0.5);
            \coordinate (O) at (0,-\cubey*0.5,-\cubez*0.5);
            
			\draw[airforceblue, thick, fill=cyan!70] (0,0,0) -- ++(-\cubex*0.5,0,0) -- ++(0,-\cubey,0) -- ++(+\cubex*0.5,0,0) -- cycle;
			\draw[airforceblue!40, thick, fill=cyan!10] (-\cubex*0.5,0,0) -- ++(-\cubex*0.5,0,0) -- ++(0,-\cubey,0) -- ++(+\cubex*0.5,0,0) -- cycle;
			\draw[airforceblue, thick, fill=red!40] (0,0,0) -- ++(0,0,-\cubez) -- ++(0,-\cubey,0) -- ++(0,0,\cubez) -- cycle;
			\draw[airforceblue, thick, fill=cyan!70] (0,0,0) -- ++(-\cubex*0.5,0,0) -- ++(0,0,-\cubez) -- ++(\cubex*0.5,0,0) -- cycle;
			\draw[airforceblue!40, thick, fill=cyan!10] (-\cubex*0.5,0,0) -- ++(-\cubex*0.5,0,0) -- ++(0,0,-\cubez) -- ++(\cubex*0.5,0,0) -- cycle;
			
			\draw[airforceblue, thick, fill=cyan!70] (-\cubex*0.5,0,0) -- ++(0,0,-\cubez) -- ++(-\cubex*0.5,0,+\cubez) -- cycle;
			\draw[airforceblue, thick, fill=cyan!70] (-\cubex*0.5,0,0) -- ++(0,-\cubey,0) -- ++(-\cubex*0.5,+\cubey,0) -- cycle;
			
			\fill[pattern=north east lines, thick, pattern color=gray] (I) -- (A) -- (J) -- (B) -- cycle;
                 
            
            \node at (E) [above] {{\scriptsize 1}};
            \node at (J) [above] {{\scriptsize 2}};
            \node at (G) [above] {{\scriptsize 3}};
            \node at (A) [above] {{\scriptsize 4}};
            \node at (H) [above] {{\scriptsize 5}};
            \node at (B) [left]  {{\scriptsize 6}};
            \node at (D) [below] {{\scriptsize 7}};
            \node at (I) [below] {{\scriptsize 8}};
            \node at (C) [below] {{\scriptsize 9}};
            \node at (L) [left] {{\scriptsize 10}};
            \node at (N) [left,xshift=-4pt,yshift=1pt] {{\scriptsize 11}};
            \node at (M) [left,xshift=-4pt,yshift=1pt] {{\scriptsize 12}};
            \node at (K) [left] {{\scriptsize 13}};
            \node at (F) [below] {{\scriptsize 14}};
            \node at (O) [right,xshift=-2pt] {{\scriptsize 15}};

			\draw[airforceblue, thick] (D) -- (E);
			\draw[airforceblue, thick] (A) -- (F);
			
			\draw[airforceblue, thick] (H) -- (E);
			\draw[airforceblue, thick] (J) -- (A);
			
			\draw[airforceblue!40, thick] (G) -- (H);
			\draw[airforceblue, thick] (H) -- (D);
			\draw[airforceblue!40, thick] (J) -- (B);
			
			\draw[airforceblue!40, thick] (B) -- (I);
			\draw[airforceblue, thick] (I) -- (A);
			\draw[airforceblue!40, thick] (C) -- (H);
			
			\draw[airforceblue, line width=0.75mm] (J) -- (H) -- (I);
			\draw[airforceblue, line width=0.75mm] (N) -- (H) -- (L);	
			
			\draw[airforceblue, line width=0.75mm] (B) -- (A);
			\draw[airforceblue, line width=0.75mm] (M) -- (H) -- (K);
			
			\draw[red!60, thick] (J) -- (B) -- (I);

			\fill[black!70] (A) circle (1pt);
			\fill[black!70] (B) circle (1pt);
			\fill[black!40] (C) circle (1pt);
			\fill[black!70] (D) circle (1pt);
			\fill[black!70] (E) circle (1pt);
			\fill[black!70] (F) circle (1pt);
			\fill[black!40] (G) circle (1pt);
			\fill[black!70] (H) circle (1pt);
			\fill[black!70] (I) circle (1pt);
			\fill[black!70] (J) circle (1pt);
			\fill[black!70] (K) circle (1pt);
			\fill[black!70] (L) circle (1pt);
			\fill[black!70] (M) circle (1pt);
			\fill[black!70] (N) circle (1pt);
			\fill[black!70] (O) circle (1pt);
			
\end{tikzpicture}
    \caption{Case 2}
    \label{fig:case2}
    \end{subfigure}
    \caption{Cancellation of the line integrals}
    \label{fig:schematic_cuboid2}
\end{figure}
Two near-field cases are considered: (1) node 13, located well within the near-field, and (2) node 5, positioned exactly at the boundary between the near-field and far-field. In the first case (see Figure \ref{fig:case1}), all the support elements of the shape function at node 13 are represented by shaded triangles. All the contributing line integrals (thickened lines) are shared between two elements and cancel out, as they are computed in opposite directions. The remaining line integrals, including those along the near-field boundary (red lines), do not contribute since the shape function vanishes there, leading to complete cancellation.

For node 5, which lies on the near-field boundary, the bounding box of the cluster is extended to fully include all its support elements within the near-field (Figure \ref{fig:case2}). As before, all contributing line integrals cancel out as they are shared between two elements. This demonstrates that line integrals are unnecessary for the near-field computation, even though the near-field itself is an open surface.

\section{Conclusions}
\label{sec6}

 The regularization of singular integrals using partial integration and Stokes' theorem is a well-established technique. However, for half-space problems, the closed surface assumption of Stokes' theorem is violated, making the line integral terms of the right-hand side necessary. In this study, the necessary line integrals for the double-layer potential and hypersingular operator were presented and their role in half-space problems was investigated. Interestingly, in the presented example, the collocation approach worked unexpectedly well, without requiring the inclusion of these new line integrals.

In contrast, the Galerkin approach, which involves the hypersingular operator, did not perform adequately for half-space problems without the inclusion of the line integrals. Only when the line integrals were incorporated did the Galerkin approach produce the correct solution, as evidenced by comparison with Boussinesq's analytical solution. 

For FMM-BEM, two modified versions of FMM were introduced to address the issue of the near-field behaving as an open surface. However, both collocation and Galerkin formulations show that the line integrals cancel in the near-field, as the bounding box is extended to fully enclose the supports of the shape functions. The results indicate that neither the line integrals nor the far-field regularization are necessary, as they contribute little to the solution. In fact, the original FMM suffices, since near-field regularization does not violate the assumptions of Stokes' theorem.

\bigskip
\noindent\textbf{Acknowledgement} \hspace*{0.5em}This work is supported by the joint DFG/FWF Collaborative Research Centre CREATOR (DFG: Project-ID 492661287/TRR 361; FWF: 10.55776/F90) at TU Darmstadt, TU Graz and JKU Linz.


\bibliography{library.bib}
\bibliographystyle{ieeetr}

\end{document}